\documentclass[fleqn]{article}

\usepackage{arxiv}

\usepackage[utf8]{inputenc} 
\usepackage[T1]{fontenc}    
\usepackage{hyperref}       
\usepackage{url}            
\usepackage{booktabs}       
\usepackage{amsfonts}       
\usepackage{nicefrac}       
\usepackage{microtype}      
\usepackage{lipsum}		
\usepackage{graphicx}
\usepackage{natbib}
\usepackage{doi}

\usepackage{float}
\usepackage{amsmath, latexsym, epsf, amssymb, framed, bbm, bm, amsthm, upref, bbold, algorithm, algorithmic}
\usepackage{xcolor}

\usepackage{makecell}

\allowdisplaybreaks

\newcommand{\Problem}[2]{
\textbf{Problem} \textsc{#1}(#2):
}
\DeclareMathOperator*{\argmax}{arg\,max}

\title{A Combined Convex Model for Travel Demand Forecasting with Hierarchical Extended Logit Model}

\author{ Youngseo Kim\\
	Civil and Environmental Engineering\\
	Cornell University\\
	\texttt{yk796@cornell.edu} \\
	\And
	Samitha Samaranayake \\
	Civil and Environmental Engineering\\
	Cornell University\\
	\texttt{samitha@cornell.edu} \\
	\AND
	Damon Wischik \\  
    Computer Science and Technology\\
    University of Cambridge\\
	\texttt{damon.wischik@cl.cam.ac.uk} \\
}

\renewcommand{\headeright}{Kim, Samaranayake, and Wischik}
\renewcommand{\undertitle}{ }
\renewcommand{\shorttitle}{ }

\hypersetup{
pdftitle={A template for the arxiv style},
pdfsubject={q-bio.NC, q-bio.QM},
pdfauthor={Youngseo Kim, Samitha Samaranayake, Damon Wischik},
pdfkeywords={First keyword, Second keyword, More},
}

\begin{document}
\maketitle

\begin{abstract}
The travel demand forecasting model plays a crucial role in evaluating large-scale infrastructure projects, such as the construction of new roads or transit lines. While combined modeling approaches have been explored as a solution to overcome the problem of input and output discrepancies in a sequential four-step modeling process, previous attempts at combined models have encountered challenges in real-world applications, primarily due to their limited behavioral richness or computational tractability. In this study, we propose a novel convex programming approach and present a key theorem demonstrating that the optimal solution is the same as the one from solving the hierarchical extended logit model. This model is specifically designed to capture correlations existing in travelers' choices, including similarities among transport modes and route overlaps. The convex property of our model ensures the existence and uniqueness of solutions and offers computational efficiency. The advantages of our proposed model are twofold. Firstly, it provides a single unifying rationale (i.e., utility maximization) that is valid across all steps. Secondly, its combined nature allows a systematic approach to handling observed data, enabling a more realistic representation of reality. By addressing the limited behavior richness of previous combined models, our convex programming-based approach shows promise in enhancing the accuracy and applicability of travel demand forecasting, thereby aiding in the planning and decision-making processes for infrastructure projects.
\end{abstract}


\section{Introduction}

Travel demand forecasting plays a crucial role in transportation planning, guiding the development of efficient and reliable transportation systems. The traditional sequential four-step model has long been the cornerstone of travel demand forecasting, encompassing trip generation, trip distribution, modal split, and traffic assignment \citep{de2011modelling}. The most commonly used models for each step are regression, gravity, logit models, and the concept of Wardrop's user equilibrium. While this model has been widely applied, critics point out its limitations due to the absence of a universal framework. Nevertheless, the four-step modeling process serves as the primary tool for forecasting future demand at the link level and assessing the performance of a transportation system. It is particularly suitable for evaluating large-scale infrastructure projects such as the construction of new roads or transit lines. 

In order to evaluate a new transportation system and estimate future demand, the four-step modeling process is typically conducted in two stages \citep{mcnally2007four}. The first stage is referred to as \textit{data assimilation}, during which the models used in each step are evaluated, calibrated, and validated using data from the existing transportation system. To give a specific example, the parameters of the regression (in trip generation), gravity (in trip distribution), and logit (in modal split) models are calibrated using the household survey data, taking into account the current transportation system. Similarly, the parameters associated with the volume-delay function (in traffic assignment) are calibrated using observed traffic data, including traffic counts and speeds. The second stage is \textit{future estimation}. In this stage, the estimated future productions and attractions are allocated to the new transportation network, which includes the newly constructed roads or transit lines. With the parameters fixed from the first stage, the same four-step procedure is then conducted to estimate the future demand and assess the performance of the updated transportation system.

However, the sequential four-step demand modeling process has an inherent weakness. The four-step model encounters challenges due to discrepancies in travel times and congestion effects at different steps of the sequential process \citep{boyce2002sequential, garrett1996transportation}. In the sequential model, an output of one step serves as an input for subsequent steps. However, the assumed values at one stage often do not align with the output values from another stage. To address this issue, equilibrated link travel times and costs from the trip assignment step are used to update the trip distribution and modal split steps in a feedback loop. However, there is no convergence guarantee. The process requires time-consuming adjustments of parameters until the results align with expectations, relying on expert judgment that may be arbitrary and ad hoc. Recognizing these issues, many researchers have made efforts to develop combined models that can simultaneously consider travelers' choices across multiple stages, resulting in more coherent outcomes. While a sequential model only allows for equilibration in the route choice process as shown in Figure \ref{fig:comparison} (a), a combined model enables equilibration in all steps as in Figure \ref{fig:comparison} (b). This better reflects the real-world scenario, as increased travel time on congested roads can impact travelers' route, mode, or even destination choices. 

Despite the advancement of four-step models, their practical application has been limited, partly due to the lack of standard software and agency inertia \citep{mcnally2007four}. Through an extensive review of the literature, we have determined that the absence of standardized software arises from the scarcity of models that can strike an effective balance between behavioral richness and computational tractability at a practical level. Models that are tractable with convexity use the simplest behavior model (i.e., multinomial logit model) \citep{oppenheim1995urban}, while models that relax undesirable assumptions lose convexity \citep{yao2014general, zhou2009}. To bridge this gap, our objective is to propose a combined model that offers a higher level of behavioral richness compared to other models that can be cast as convex programs.

The research makes a substantial contribution in three main areas. Firstly, we introduce a hierarchical extended logit model that effectively captures the sequential decision-making process of destination, mode, and route choice. By incorporating the nested logit model within the mode choice stage, we account for the correlation among modes. Additionally, the inclusion of the path-size logit model in the route choice step enables us to consider route overlapping, resulting in a desirable level of behavior richness. Secondly, the proposed model maintains the advantageous convexity of the formulation while preserving the desirable behavior richness. The convexity property ensures the scalability of the model to real-world networks, making it practical for application in various scenarios. Thirdly, we propose a two-stage travel demand forecasting framework within the combined model: data assimilation and future estimation. Notably, we suggest a formulation for the data assimilation stage that suits a combined approach. The most probable parameters are obtained from the dual variables of the constraints for available observations, instead of maximum likelihood estimation. Conducting data assimilation using the combined model offers multiple benefits, including more reliable and coherent outcomes by leveraging diverse sources of data for each of the four steps. Importantly, to the best of our knowledge, this is the first comprehensive model to properly incorporate real-world data into combined travel demand modeling, distinguishing it from previous approaches.

\subsection{Sequential four step model}
\begin{figure}[H]
    \centering
    \includegraphics[width=16cm]{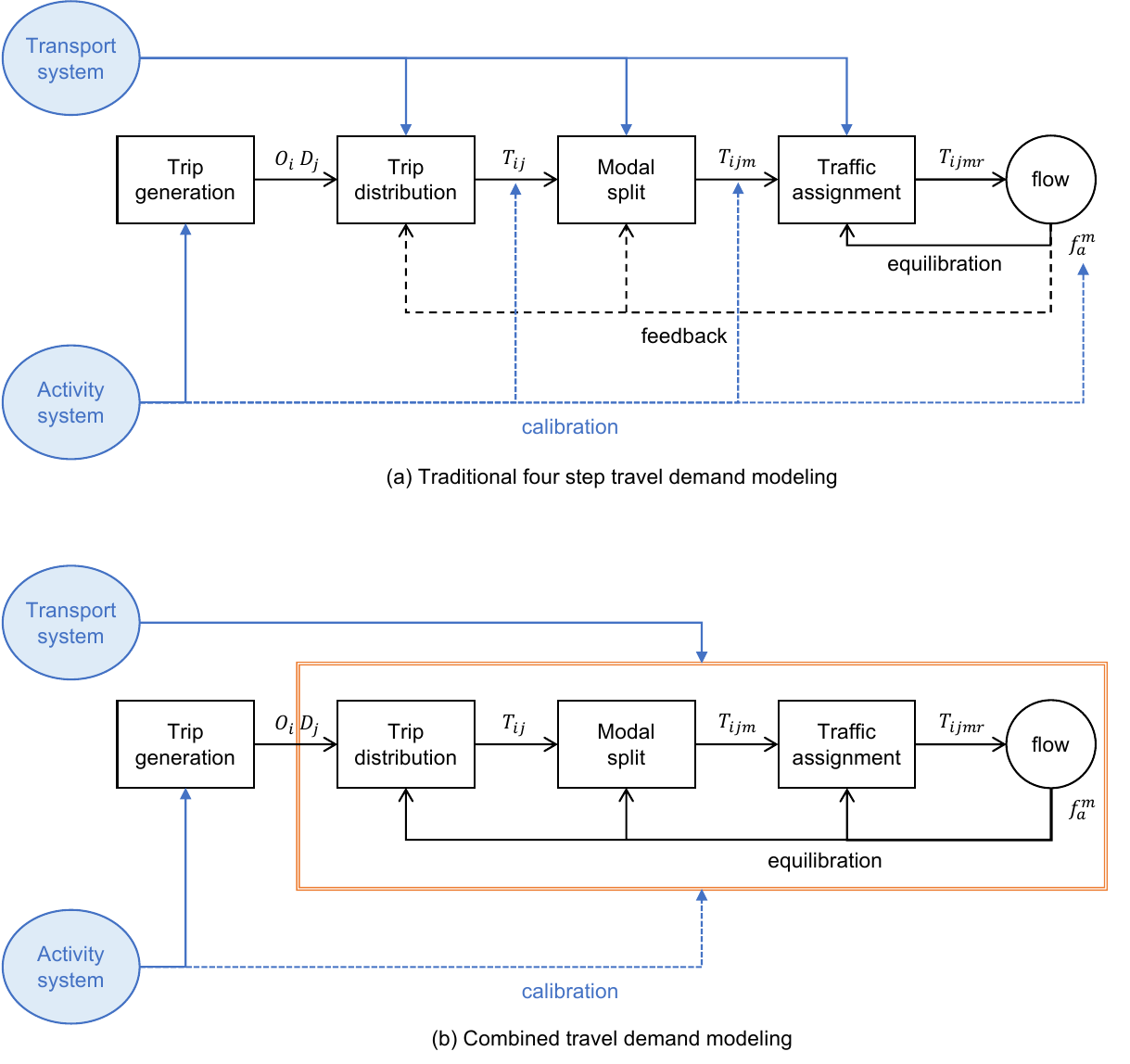}
    \caption{Traditional four step model (modified from \citet{mcnally2007four}) \textit{vs} our framework}
    \label{fig:comparison}
\end{figure}

Figure \ref{fig:comparison} (a) provides a concise overview of the traditional four-step demand modeling. Within this framework, two exogenous factors play a crucial role: the transportation system and the activity system. The transportation system is represented by network graphs consisting of links and nodes with associated attributes such as length, speed, and capacity. On the other hand, the activity system is typically characterized by socioeconomic, demographic, and land use data that are defined for specific spatial units, such as Traffic Analysis Zones (TAZs) in the United States. To capture the activity system, household survey data is collected as a representative sample. 

We now examine each of the four steps individually, keeping in mind that the explanations provided represent common modeling practices in the United States. To maintain consistency with the notations, we provide the corresponding variables in our model in parentheses, which will be defined in subsequent sections. The first step is trip generation. Utilizing household information (e.g., car ownership, income) and zonal data (e.g., employment, land use), we estimate separate productions ($O_i$) and attractions ($D_j$) using regression-based or category models. These productions and attractions can be generated for specific times of the day and different trip purposes. The second step is trip distribution. With the given productions and attractions, an origin-destination (OD) trip matrix ($T_{ij}$) is completed using the Gravity model, which reflects the underlying travel impedance (travel time and cost). The third step involves mode choice. Using logit models, OD trip matrices for each mode ($T_{ijm}$) are completed to represent the relative proportions of trips made by each travel mode. The fourth step is route choice. Modal OD trip matrices are assigned to mode-specific networks. In this step, our objective is to attain a user equilibrium where all paths utilized for a given OD pair have equal impedance. The basic user equilibrium solution is obtained using the Frank-Wolfe algorithm, which involves computing minimum paths, performing all-or-nothing assignments to these paths, and obtaining flow for each link.

The paper is organized as follows. In the preliminary and literature review section, we discuss equivalent optimization models for trip distribution, modal split, and trip assignment. Specifically, we explore how the gravity model and logit model can be reformulated as entropy maximization frameworks. Additionally, we review existing combined models and highlight the distinctive features of our proposed approach. In the methodology and results section, we delve into our main contribution by introducing entropy maximization formulations for the nested logit model and hierarchical multinomial logit model. We also propose a novel two-stage formulation for the hierarchical extended logit model. Lastly, the conclusion section provides a concise summary of our findings, emphasizing their significance. We discuss the implications of our research and suggest future avenues for advancing travel demand forecasting and transportation planning.

\section{Preliminary and literature review}\label{section:preliminary}

\newtheorem{theorem}{Theorem}
\newtheorem{assumption}{Assumption}
\newtheorem{corollary}{Corollary}
\newtheorem{proposition}{Proposition}
\newtheorem{lemma}{Lemma}
\newtheorem{definition}{definition}

\newcommand{\tauvec}[0]{\ensuremath{\boldsymbol{\tau}}}
\newcommand{\betavec}[0]{\ensuremath{\boldsymbol{\beta}}}
\newcommand{\alphavec}[0]{\ensuremath{\boldsymbol{\alpha}}}
\newcommand{\gammavec}[0]{\ensuremath{\boldsymbol{\gamma}}}
\newcommand{\varepsilonvec}[0]{\ensuremath{\boldsymbol{\varepsilon}}}
\newcommand{\nest}[0]{\ensuremath{M}}
\newcommand{\nests}[0]{\ensuremath{\mathcal{M}}}

The term \textit{entropy}, which quantifies the amount of energy in a thermodynamic system, was first introduced by \citet{shannon1948mathematical} in the field of information theory. Shannon's concept of entropy maximization is to obtain a probability distribution that complies with observed information regarding an unknown distribution. In this section, we focus on the theory of entropy maximization developed in transportation planning. For a more comprehensive understanding of entropy maximization and its diverse range of applications, readers are encouraged to refer to \citet{fang1997entropy}.

\subsection{Entropy maximization for trip distribution}

\citet{wilson1969use} propose an entropy maximization model as an alternative to the conventional gravity model for the trip distribution step. The entropy maximization approach aims to determine the most likely distribution of trips and serves as a theoretical foundation for the gravity model. Subsequently, the author proposes generalized distribution models that incorporate multiple transportation modes, thereby uncovering an implicit modal split model.

\subsubsection{Gravity model for a single mode}

We first focus on the gravity model for a single mode of transport. We extend it to several transport modes later. The region is divided into zones: origin zone $i \in I$ and destination zone $j \in J$. The gravity model is used to estimate $T_{ij}$, the number of trips between $i$ and $j$. Here, we deal with aggregated decision makers at spatial zone levels, and do not distinguish individual travelers.

\begin{subequations}\label{eqn:gravity}
\begin{align}
    & T_{ij} 
    = A_i B_j O_i D_j f(c_{ij}), & \forall i,j \\
    & \text{where } & \nonumber \\
    & A_i = \frac{1}{\sum_{j} B_j D_j f(c_{ij})} , & \forall i \\
    & B_j = \frac{1}{\sum_{j} A_i O_i f(c_{ij})}, & \forall j
\end{align}
\end{subequations}

\noindent where $O_i$ be the total number of trip origins at $i$ (i.e., trip production). $D_j$ be the total number of trip destinations at $j$ (i.e., trip attraction). $c_{ij}$ be the generalized cost (e.g. a linear sum of fares and travel times). $f(\cdot)$ be a decreasing cost function.

\subsubsection{Entropy maximization model for a single mode}
Now, we suggest the entropy maximization model to get the most probable distribution of trips. We eventually show that the gravity model can be obtained by solving the entropy maximization model. The basic assumption is that the probability of $[T_{ij}]_{i \in I, j \in J}$ is proportional to the number of states of the system with such distribution. Based on combinatorial theory, we express the number of distinct arrangements, $w([T_{ij}]_{i \in I, j \in J})$.

\begin{equation}
    w([T_{ij}]_{i \in I, j \in J}) = \frac{T!}{\Pi_{ij} T_{ij}!}
\end{equation}

\noindent where $T$ is the total number of trips.

We find a set of $T_{ij}$ which maximize $\log w([T_{ij}]_{i \in I, j \in J})$ in order to get the most probable distribution of trips.

\Problem{{\fontfamily{cmr}\selectfont MostProbable}}{$[O_i]_{i \in I}, [D_j]_{j \in J}, C$}

\begin{subequations}
\begin{align}
\max_{[T_{ij}]_{i \in I, j \in J}} \quad &\log w([T_{ij}]_{i \in I, j \in J}) \label{mostprob:logw} \\
    \text{s.t.} \quad & \sum_{j} T_{ij} = O_i, \quad \forall i \quad [\lambda_i] \label{mostprob:sum_origin}\\
    & \sum_{i} T_{ij} = D_j, \quad \forall j \quad [\mu_j] \label{mostprob:sum_destin}\\
    & \sum_{i} \sum_{j} T_{ij} c_{ij} = C \quad [\beta] \label{mostprob:sum_cost}
\end{align}
\end{subequations}

\noindent Objective function in Equation (\ref{mostprob:logw}) finds the trip distribution with the maximum probability. Using Stirling's approximation (i.e., $\log T! = T \log T - T$), the objective function can be written as 
\begin{align}
    \quad & \log w([T_{ij}]_{i \in I, j \in J}]_{i \in I, j \in J}) = \log T! - \sum_{ij} \log T_{ij}! \approx \log T! - \sum_{i}\sum_{j} T_{ij} \log T_{ij} + \sum_{i}\sum_{j} T_{ij}
\end{align}
\noindent Constraint (\ref{mostprob:sum_origin}) and (\ref{mostprob:sum_destin}) ensure that the total number of trip from origin $i$ and to destination $j$ should be equal to given observations. Constraint (\ref{mostprob:sum_cost}) ensures to satisfy the budget constraint where $C$ is the total amount spent on these trips in the region. Dual variables are listed in square brackets.

The Lagrangian function $\mathcal{L}$ is
\begin{align}
    {\mathcal{L}} =& - \sum_{i}\sum_{j} T_{ij} \log T_{ij} + \sum_{i}\sum_{j} T_{ij} \nonumber \\ & + \sum_{i}{\lambda_i} (O_i - \sum_{j} T_{ij}) 
    + \sum_{j}{\mu_j} (D_j - \sum_{i} T_{ij}) 
    + \beta (C - \sum_{i} \sum_{j} T_{ij})
\end{align}

\noindent Note that we call this approach an entropy maximization method due to the term $- \sum_{i}\sum_{j} T_{ij} \log T_{ij}$. Using $p_{ij} = T_{ij}/T$, we obtain $H = - \sum_{i} \sum_{j} p_{ij} \log p_{ij}$, which is the information theorist's definition of entropy.

By solving the first order condition, $\frac{\partial {\mathcal{L}}}{\partial T_{ij}} = 0, \forall i, j$, we recover the gravity model. 

\begin{subequations}
\begin{align}
    & \widehat{T_{ij}} = A_i B_j O_i D_j \exp(-\beta c_{ij}), & \forall i,j \\
    \text{where } & A_i = \frac{1}{\sum_{j} B_j D_j \exp(-\beta c_{ij})} \\
                  & B_j = \frac{1}{\sum_{j} A_i O_i \exp(-\beta c_{ij})}
\end{align}
\end{subequations}

\noindent Thus, the most probable distribution of trips is the same as the one normally recognized as the gravity model distribution suggested in Equation (\ref{eqn:gravity}), where a decreasing cost function $f(c_{ij})$ is set to $\exp(-\beta c_{ij})$. This statistical derivation constitutes a new theoretical base for the gravity model. Note that $C$ does not need to be known as $\beta$ is given in practice. Also, if there were no constraints for observation (i.e., without Constraint (\ref{mostprob:sum_origin}) and (\ref{mostprob:sum_destin})), $T_{ij}$ would have an equal share of the total number of trips. In other words, it is the constraints that make this trip distribution not trivial. 

\subsubsection{Entropy maximization model for several modes}
We introduce several transport modes with notations extended to multiple modes. $m$ denotes a mode and $m \in \mathcal{M}$ where $\mathcal{M}$ is the set of available modes. $T^{m}_{ij}$ denotes the number of trips between $i$ and $j$ by mode $m$. $c^m_{ij}$ denotes the generalized cost of traveling from $i$ to $j$ by mode $m$.

Following the same procedure as before, we get the following equations. 

\begin{align}
    & \widehat{T^m_{ij}} = A_i B_j O_i D_j \exp(-\beta c^m_{ij}), & \forall i,j,m \\
    \text{where } & A_i = \frac{1}{\sum_{j} B_j D_j \exp(-\beta c^m_{ij})} \\
                  & B_j = \frac{1}{\sum_{j} A_i O_i \exp(-\beta c^m_{ij})}
\end{align}

Accordingly, the modal split can be written as follows. 

\begin{align}\label{eqn:entropy_modalsplit}
    \frac{T^{m}_{ij}}{\sum_{m' \in \mathcal{M}}T^{m'}_{ij}} = \frac{\exp(-\beta c^m_{ij})}{\sum_{m' \in \mathcal{M}} \exp(-\beta c^{m'}_{ij})}, & \forall i,j,m
\end{align}

We have shown that the entropy maximization model for several modes indeed embeds a logit expression of modal split. In the following section, we investigate more deeply how behavior demand models and entropy maximizing models are related. 

\subsection{Entropy maximization for modal split}

We explore entropy maximizing models that substitute logit choice models. Using the family of logit models, especially the multinomial logit model (MNL), are the standard way to model traveler mode choices. We show that MNL can be derived from the entropy maximization model. 

\subsubsection{Multinomial logit model and maximum likelihood estimation}

\citet{mcfadden1973conditional} was the first to formulate travel mode and location decisions as problems in micro-economic consumer choices. In the random utility theory, the utility of an alternative $U^h_{m}$ can be expressed as the sum of the systematic component $V^h_m$ and the error component $\varepsilon^h_{m}$ associated with joint random variation across both individuals $h$ and alternatives $m$. $X^h_{mk}$ is the $k$-th trip related attribute (e.g. travel cost or travel time) for alternative $m$, and $\beta_{k}$ is a parameter for $k$-th attribute. $ASC_{m}$ is the alternative specific constant for alternative $m$. We focus on specific $ij$-pair and omit the subscript.

\begin{align}
    U^h_{m} = V^h_{m} + E^h_{m} = ASC_{m} + \sum_{k \in K} \beta_{k} X^h_{mk} + E^h_m 
\end{align}

If $\varepsilon^h_m$ follows the type I extreme value distribution (also known as Gumbel distribution) which is independent and identically distributed, the cumulative distribution function of an error term is

\begin{align}
    F([\varepsilon^h_m]_{m \in M}) = \mathbb{P} (E^h_m \leq \varepsilon^h_m, \forall m \in M) = e^{-e^{-{\theta \varepsilon^h_m}}}, \theta >0 \label{eqn:gumbel_mnl}
\end{align}

\noindent With the assumption of the Gumbel distribution, \citet{ben1985discrete} suggested that the choice probability can be expressed as a closed form logit function.

\begin{subequations}\label{eqn:mnl_logit}
\begin{align}
    P^h_{m}(\betavec; \theta) 
    &= \mathbb{P}(U^h_{m} > U^{h}_{m'}, \forall m' \neq m) \\
    & = \frac{e^{\theta V^h_{m}}}{\sum_{m'} e^{\theta V^{h}_{m'}}} 
    = \frac{e^{\theta (ASC_{m} + \sum_{k} \beta_{k} X^h_{mk})}}{\sum_{m'} e^{\theta (ASC_{m'} \sum_{k} \beta_{k} X^h_{m'k})}}, \quad \quad \forall m,h
\end{align}
\end{subequations}

\noindent With the given observation of $y^h_m$, the indicator of whether individual $h$ chooses alternative $m$, we maximize the log-likelihood function to estimate the coefficient $\betavec$.

\begin{align}
    & \bar{\betavec} 
    = \argmax_{\betavec} \log \mathcal{L} \\
    & \mathcal{L} = \sum_h \sum_m y^h_m \log P^h_m(\betavec)
\end{align}

\noindent With the estimated $\bar{\betavec}$, we obtain $\overline{V^h_{m}}$, and then obtain $p^h_{m}$ from the equation below. 

\begin{align}\label{eqn:prob_fixed_point}
    p^h_{m} &
    = \frac{e^{\theta \overline{V^h_{m}}}}{\sum_{m'} e^{\theta \overline{V^{h}_{m'}}}}, \quad \quad \forall m,h
\end{align}

\subsubsection{Satisfaction maximization model}\label{section:utility_maximize_MNL}

Satisfaction function $S(\mathbb{V})$ is defined as the expected utility that a traveler received from the set of alternatives. We assume homogenous individuals. 

\begin{align}
    S(\mathbb{V}) = \mathbb{E}[\max_{\forall m \in \mathcal{M}} \{ V_m + E_m \}] = \int_{-\infty}^{\infty} \max_{\forall m \in \mathcal{M}}( V_m + \varepsilon_m) F(\varepsilon_m) \mathrm{d} \varepsilon_m
\end{align}

\noindent By assuming that $F$ is the extreme value distribution, we have the following satisfaction function and the choice probabilities. The probabilistic formulation exhibits a multinomial logit structure as in Equation (\ref{eqn:mnl_logit}) that has an obvious microeconomic interpretation. For any positive $\theta$,

\begin{align}
    S(\boldsymbol{V}) = \frac{1}{\theta} \log \sum_{m}e^{\theta V_m}
\end{align}

\begin{align}
    p_m(V) = \frac{e^{\theta V_m}}{\sum_{m'}e^{\theta V_{m'}}} .
\end{align}

\citet{miyagi1996direct} suggested the fundamental equation called \textit{choice equation} as below.

\begin{align}\label{eqn:choice_equation}
    S(V) = \sum_{m \in \mathcal{M}} V_m p_m + H(p_m) = \sum_{m \in \mathcal{M}} V_m p_m - \frac{1}{\theta} \sum_{m \in \mathcal{M}} p_m \ln (p_m)
\end{align}

\noindent Note that $H(p_m)$ has a form of entropy function, which is a well-known measure of uncertainty. In this case, the entropy term is a direct measure of the increase in utility arising from the increase in the variety of choices \citep{miyagi1996direct}. Accordingly, we refer to $H(p_m)$ as the freedom of choice. In addition, since the entropy function appears in the satisfaction function, satisfaction maximization can be viewed as entropy maximization.

One interesting property is that the probability that maximizes the satisfaction function is the choice probability obtained by solving the fixed point problem using the MNL. The equivalent convex optimization formulation that can substitute the MNL is as follows. Suppose the deterministic utility $\overline{V_m}$ is given.

\Problem{{\fontfamily{cmr}\selectfont MaxSatisMNL}}{$[\overline{V_m}]_{m \in \mathcal{M}}, \theta$}

\begin{subequations}
\begin{align}
    \max_{[p_m]_{m \in \mathcal{M}}} & \sum_{m \in \mathcal{M}} \overline{V_{m}} p_m - \frac{1}{\theta} \sum_{m \in \mathcal{M}} p_m \ln (p_m) \\
    \text{s.t.} & \sum_{m \in \mathcal{M}} p_m = 1, 
\end{align}
\end{subequations}

\noindent By solving the first order condition, we can recover the choice probability.

\begin{align}
    \widehat{p_m} &= \frac{e^{\theta \overline{V_{m}}}}{\sum_{m'} e^{\theta \overline{V_{m'}}}}, \forall m \label{eqn:prob_max_entropy}
\end{align}

\noindent With the assumption of homogeneous individual, the choice probability in Equation (\ref{eqn:prob_max_entropy}) is the same as Equation (\ref{eqn:prob_fixed_point}) which is derived from the MNL. Note that this is a generalization of Equation (\ref{eqn:entropy_modalsplit}), expressing generalized cost ($c^m_{ij}$) more precisely with travel attributes ($\sum_{k} \overline{\beta_k}X_{mk}$).

\subsubsection{Theoretical background}

We show the theoretical background of the reason why the entropy term appears in the satisfaction maximization model. \citet{anas1983discrete} showed that the behavior demand model \citep{mcfadden1973conditional} and the entropy maximizing model \citep{wilson1969use}, should be seen as two equivalent views of the same problem. Consider the following entropy maximizing problem. 

\Problem{{\fontfamily{cmr}\selectfont MaxEntropy}}{$[\overline{y^h_{m}}]_{m \in \mathcal{M}, h \in H}, \theta$}

\begin{subequations}
\begin{align}
    \max_{[p^h_{m}]_{\forall h \in H, m \in M}} &  - \frac{1}{\theta}\sum_{h} \sum_{m} p^h_{m} \log p^h_{m} \label{maxentropy:objective}\\
    \text{s.t.} & \sum_{m} p^h_{m} = 1; \quad \forall h \quad [\lambda_h] \label{maxentropy:prob_sum}\\
    & \sum_{h} p^h_{m} = \sum_{h} \overline{y^h_{m}}; \quad \forall m \quad [\gamma_m] \label{maxentropy:gamma_m}\\
    & \sum_{m} \sum_{h} p^h_{m} X^h_{mk} = \sum_{m} \sum_{h} \overline{y^h_{m}} X^h_{mk}, \quad \forall k \quad [\alpha_k] \label{maxentropy:alpha_0}
\end{align}
\end{subequations}

\noindent Objective function (\ref{maxentropy:objective}) seeks the most random (entropy maximizing) predictions of $p^h_m$. Constraint (\ref{maxentropy:prob_sum}) requires the choice probability of choosing each alternative should be summed up to one for each individual. Constraint (\ref{maxentropy:gamma_m}) ensures that the predicted number of choosing each alternative should equal to the actual number of choosing it. Constraint (\ref{maxentropy:alpha_0}) ensures that the expectation of the aggregated deterministic value of each alternative should equal the observed aggregated value. Constraint (\ref{maxentropy:gamma_m}) and (\ref{maxentropy:alpha_0}) both ensure that these predictions replicate the aggregated observations on the entire system.

Similar to what we did in the previous section, by forming the Lagrangian function of the above formulation and by solving the first order condition, we obtain the form of MNL.  

\begin{align}
    \widehat{p^h_{m}}  
    = \frac{e^{\theta (\gamma_m + \sum_{k}\alpha_k X^{h}_{mk})}}{\sum_{m'} e^{\theta (\gamma_{m'} + \sum_{k}\alpha_k X^{h}_{m'k})}} 
\end{align}

\noindent \citet{anas1983discrete} proved that the MNL model can be identically estimated via maximum likelihood that finds $\betavec$ and $\textbf{ASC}$ in Equation (\ref{eqn:mnl_logit}), or via solving {\fontfamily{cmr}\selectfont MaxEntropy} to find Lagrangian multipliers $\alphavec$ and $\gammavec$, thus $\betavec = \alphavec$ and $\textbf{ASC} = \gammavec$. The maximum entropy estimators reproduce rational user behavior. 

In this section, we have investigated the equivalence between the entropy maximization model and the gravity model employed in the trip distribution step. Furthermore, we examine the equivalence between the satisfaction maximization model and the MNL in the modal split step. Additionally, it is widely acknowledged that \citet{beckmann1956studies} proposed a convex optimization formulation for the Wardropian user equilibrium assignment problem in the traffic assignment step. These findings highlight the existence of convex optimization problems that can replace the traditional fixed point approaches in the trip distribution, modal split, and traffic assignment steps. 

\subsection{Combined modeling approach}

\begin{table}[ht]
\caption{Literature for Combined Models for Four Step Travel Demand Modeling}
\centering
\label{tab:literature}
\begin{tabular}{@{}lcccc@{}}
\toprule
\textbf{}                                        & \textbf{\begin{tabular}[t]{@{}c@{}}Trip \\ generation \end{tabular}} & \textbf{\begin{tabular}[t]{@{}c@{}}Trip  \\ distribution\\ (Gravity \\ model) \end{tabular}} & \textbf{\begin{tabular}[t]{@{}c@{}}Modal \\ split\\ (Multinomial \\ logit model)\end{tabular}} & \textbf{\begin{tabular}[t]{@{}c@{}}Traffic \\ assignment\\ (Wardrop \\ user equilibrium)\end{tabular}} \\ \midrule
\textbf{Framework}                 & \textbf{Utility}                & \textbf{Entropy}                                                                     & \textbf{Satisfaction}                                                                       & \textbf{Beckmann}                                                                                   \\ \midrule
\citet{wilson1969use}            &                          & O                                                                                    &                                                                                             &                                                                                                     \\
\citet{anas1983discrete}         &                          &                                                                                      & O                                                                                           &                                                                                                     \\
\citet{beckmann1956studies}      &                          &                                                                                      &                                                                                             & O                                                                                                   \\
\midrule
\citet{evans1976derivation}      &                          & O                                                                                    &                                                                                             & O                                                                                                   \\
\citet{florian1975combined}      &                          & O                                                                                    &                                                                                             & O                                                                                                   \\
\citet{oppenheim1993equilibrium} &                          & O                                                                                    &                                                                                             & O                                                                                                   \\
\citet{florian1977traffic}       &                          &                                                                                      & O                                                                                           & O                                                                                                   \\
\citet{abdulaal1979methods}      &                          &                                                                                      & O                                                                                           & O                                                                                                   \\
\citet{fernandez1994network}     &                          &                                                                                      & O                                                                                           & O                                                                                                   \\
\citet{garcia2005network}        &                          &                                                                                      & O                                                                                           & O                                                                                                   \\
\citet{florian1978}              &                          & O                                                                                    & O                                                                                           & O                                                                                                   \\
\citet{friesz1981equivalent}     &                          & O                                                                                    & O                                                                                           & O                                                                                                   \\
\citet{safwat1978}               & O                        & O                                                                                    & O                                                                                           & O                                                                                                   \\ \midrule
\textbf{Framework}                               & \multicolumn{4}{c}{\textbf{Utility maximization}}                                                                                                                                                                                                                                                                   \\ \midrule
\citet{oppenheim1995urban}       & O                        & O                                                                                    & O                                                                                           & O                                                                                                   \\
\citet{yao2014general}           & O                        & O                                                                                    &                                                                                             & O                                                                                                   \\
\citet{zhou2009}                 & O                        & O                                                                                    & O                                                                                           & O                                                                                                   \\
\textbf{Ours}                                    & \textbf{O}               & \textbf{O}                                                                           & \textbf{O}                                                                                  & \textbf{O}                                                                                          \\ \bottomrule
\end{tabular}
\end{table}


As in Table \ref{tab:literature}, building upon the influential studies by \citet{wilson1969use}, \citet{anas1983discrete}, and \citet{beckmann1956studies} discussed in the previous section, subsequent literature has integrated relevant terms (i.e., entropy, satisfaction, Beckmann function) into objective functions. These terms are selected based on the specific step that researchers aim to combine. The previous literature can be categorized into four groups: (a) combining trip distribution and traffic assignment \citep{evans1976derivation, florian1975combined, oppenheim1993equilibrium}, (b) combining modal split and traffic assignment \citep{florian1977traffic, abdulaal1979methods, fernandez1994network, garcia2005network}, (c) combining trip distribution, modal split, and traffic assignment \citep{florian1978, friesz1981equivalent}, and (d) combining trip generation, distribution, modal split, and traffic assignment \citep{safwat1978}.

Another stream of research has tried to express all steps in the \textit{utility maximization} framework \citep{oppenheim1995urban, zhou2009, yao2014general}. Note that the term entropy maximization and utility maximization are interchangeable in their setting. Both terms emphasize that an equivalent optimization problem is solved instead of the fixed point problem. In the utility maximization framework, an individual traveler acts as a rational consumer of urban trips, and their choices are governed by the random utility theory. The solution to the problem corresponds to collective utility maximization as well as individual utility. The utility maximization framework naturally incorporates the trip generation step because the framework assumes trip demand at the origin as given values and builds a model for destination choice. This approach is reasonable considering that the production model, in general, is more reliable, and the attractions are often normalized by the number of productions \citep{mcnally2007four}. The way that \citet{safwat1978} integrates the trip generation step is the utility maximization framework with an accessibility measure.

\citet{oppenheim1995urban} is the first to propose a convex optimization representing a hierarchical structure of the travel-destination-mode-route choice and simultaneous decision making process. For each choice step, the MNL is used. Subsequent researchers built upon this research \citep{zhou2009, yao2014general}. \citet{zhou2009} presented an unconstrained optimization formulation that embeds a hierarchical structure, providing the flexibility to handle general probabilistic distributions. To tackle the challenge of non-separable link cost functions (where the link travel time of each mode depends on the flows of other modes on that link), they proposed a variational inequality formulation. While this formulation guarantees the existence and uniqueness of a solution, it necessitates iterative methods that are not computationally efficient. The authors suggested a projection algorithm specifically for solving a hierarchical logit model. However, it was acknowledged that this approach would increase computational overhead when dealing with large-scale networks. As such, the proposed method might not be as practical for real-world applications with significant network sizes.

\citet{yao2014general} aimed to address the limitations associated with the independence of irrelevant alternative (IIA) assumption inherited from the MNL. To achieve this, they proposed the spatially correlated logit model to account for unobserved similarities among destinations in the destination choice, and the path-size logit model to consider route overlapping in the route choice. However, while the extended logit model allows to capture certain behavior complexities, the formulation does not maintain convexity. As a result, only the simplest MNL was suggested as a special case that adheres to the convexity property. More importantly, their hierarchical logit model for the combined model excludes the modal split step and focuses exclusively on the road network. This specialization may limit its applicability when considering other transportation modes.

\section{Methodology and Results}

In this section, we begin by introducing the nested logit model and the hierarchical multinomial logit model within the entropy/utility maximization framework as a preliminary step. Subsequently, we put forward the hierarchical extended logit model as our main contribution. We present the two-stage hierarchical extended logit model in the entropy/utility maximization framework. This novel approach incorporates the principles of \textit{data assimilation} and \textit{future estimation} from the traditional sequential model (see Figure \ref{fig:comparison}), but with significant improvements.

\subsection{Utility maximization for nested logit model}

\subsubsection{Nested logit model}

The multinomial logit model (MNL) is a widely used choice model that assumes the independence of irrelevant alternatives (IIA), meaning the introduction or removal of an alternative should not impact the relative preference between other alternatives. Because of this assumption, MNL has a significant limitation; it often fails to capture realistic substitution patterns among alternatives. This is evident in scenarios like the well-known \textit{blue bus and red bus} paradox. Imagine a transportation market with private cars and red buses, each with a 50\% market share. Now, a new alternative is introduced - a blue bus. According to the IIA assumption, the blue bus should take share equally from both cars and red buses. For example, the market share for each mode can be 33.3\%. However, this prediction seems unreasonable because the introduction of the new bus disproportionately affects the market share of existing alternatives. Common sense suggests that the blue bus would mostly take share from the red bus rather than private cars. Therefore, it is more reasonable to expect that the total bus share would remain close to 50\%. 

To address this limitation and better reflect the complexities of decision-making, researchers often resort to introducing the nested logit model (NL). The NL model relaxes the IIA assumption by allowing for nested groups of alternatives, enabling a more flexible representation of substitution patterns. By incorporating the NL framework, analysts can more accurately model consumer choices in situations with multiple alternatives and better understand how the introduction of new options influences existing choices.

\begin{figure}[H]
    \centering
    \includegraphics[width=\textwidth]{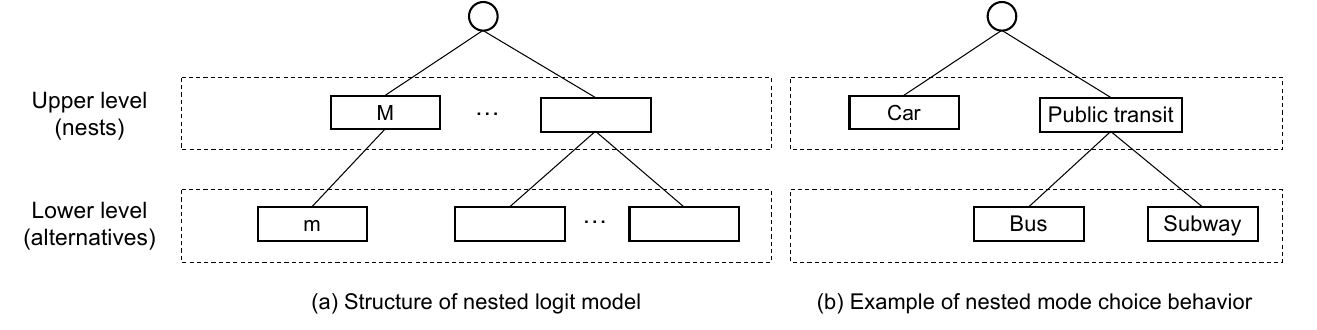}
    \caption{Nested logit model}
    \label{fig:nestedlogit}
\end{figure}

The NL model is a generalization of the MNL model that accounts for the correlation among alternatives within nested groups \citep{wen2001generalized}. We consider the 2-layer NL model as in Figure \ref{fig:nestedlogit} for simplicity, but this can be extended to $n$-layer NL models. The main difference from the MNL is the existence of nests that denote similar alternatives within a group. Let $m$ be an alternative, $M$ be a nest, and $\mathcal{M}$ be a set of all alternatives.  The upper layer contains nest $\nest \in \nests$. Function $B(m)$ links between the upper layer and the lower layer and indicates the nest for $m$. The NL model introduces the concept of dissimilarity among alternatives within a nest. $\tau_{\nest}$ denotes the dissimilarity factor for the nest $\nest$, and it is equal to $\sqrt{1-p_{\nest}}$ where $p_{\nest}$ is the coefficient of correlation, i.e., measure or similarity. When there is no correlation (i.e., $p_{\nest} = 0$), the dissimilarity factor has the highest value (i.e., $\tau_{\nest} = 1$), which is in the case that the NL model equals the MNL model. 

The random utility of the NL model can be written as follows. 

\begin{subequations}
\begin{align}
    & U_{m} = V_{m} + E_{m} \\
    & F([\varepsilon_m]_{m \in \mathcal{M}}; \theta, [\tau_{M}]_{M \in \mathcal{M}}) = e^{-\sum_{M \in \mathcal{M}} [\sum_{m \in M} e^{-\theta \varepsilon_m/{\tau_{M}}}]^{\tau_{M}}} \label{eqn:NL_cdf}
\end{align}
\end{subequations}

\noindent Equation (\ref{eqn:NL_cdf}) is the cumulative distribution function embedding the correlation within nests.  

The choice probability can be expressed as below. 

\begin{subequations}
    \begin{align}
        & p_M = \frac{e^{IV_{M}}}{\sum_{\nest \in \nests}e^{IV_{\nest}}} \label{nl:prob_nest}\\
        & p_{m/M} =\frac{e^{\theta V_{m}/{\tau_{M}}}}{\sum_{m' \in M}e^{\theta V_{m'}/\tau_{M}}} \label{nl:prob_alter}\\
        & \text{where }IV_{\nest} = \tau_M \ln \sum_{m' \in \nest }e^{\theta V_{m'}/\tau_{\nest}} \label{nl:IV} 
    \end{align}    
\end{subequations}

\begin{equation}\label{eqn:nl_2}
    \therefore p_{m} = p_{B(m)} p_{m/B(m)}
    = \frac{e^{\theta V_{m}/{\tau_{B(m)}}}}{e^{IV_{B(m)}/\tau_{B(m)}}} \times \frac{e^{IV_{B(m)} }}{\sum_{\nest \in \nests}e^{IV_{\nest}}}
\end{equation}

The probability of choosing nest $M$ is expressed as a logit model with the aggregated utility of the nest, $IV_{M}$ as in Equation (\ref{nl:prob_nest}). The inclusive value ($IV_{\nest}$) is defined in Equation (\ref{nl:IV}). The inclusive value represents the attractiveness of nest $\nest$ since it measures the summation of normalized utilities of all alternatives in nest $\nest$, which allows us to compare alternatives across nests. Equation (\ref{nl:prob_alter}) shows the probability of choosing alternative $m$ within its nest, which is expressed as a logit model with the utility term scaled by $\tau_M$. As can be seen in Equation (\ref{eqn:nl_2}), the probability of choosing alternative $m$, $p_m$, is the product of the probability to choose nest $B(m)$, $p_{B(m)}$, and the conditional probability to choose alternative $m$ in the nest $B(m)$, $p_{m/B(m)}$.   


\subsubsection{Nested logit model in entropy maximization framework}

Recall the satisfaction function for the MNL model defined in the previous section. Now, we can extend \textit{choice equation} in Equation (\ref{eqn:choice_equation}) to the NL model with alternative $m \in \mathcal{M}$.

\begin{align}\label{eqn:choice_equation_nl}
    S(V) & = \sum_{m \in \mathcal{M}} V_m p_m + H(p_m) \\
    & = \sum_{m \in \mathcal{M}} V_m p_m - \frac{1}{\theta} \sum_{m \in \mathcal{M}} \left[\tau_{B(m)} p_m \ln (p_m) + (1-\tau_{B(m)}) p_m \ln (\sum_{m' \in B(m)}p_{m'})\right]
\end{align}

\Problem{{\fontfamily{cmr}\selectfont MaxSatisNL}}{$[\overline{V_m}]_{m \in \mathcal{M}}, \theta, [\tau_{B(m)}]_{m \in \mathcal{M}}$}

\begin{subequations}\label{eqn:max_choice_equation_nl}
\begin{align}
    \max_{[p_m]_{m \in \mathcal{M}}} & \sum_{m \in \mathcal{M}} \overline{V_m} p_m - \frac{1}{\theta} \sum_{m \in \mathcal{M}} \left[\tau_{B(m)} p_m \ln (p_m) + (1-\tau_{B(m)}) p_m \ln (\sum_{m' \in B(m)}p_{m'})\right] \\
    \text{s.t.} & \sum_{m \in \mathcal{M}} p_m = 1, 
\end{align}
\end{subequations}

\begin{theorem}\label{theorem:choice_equation_nl}
    The optimal solution $[p_m]_{m \in \mathcal{M}}$ of the convex maximization problem {\fontfamily{cmr}\selectfont MaxSatisNL} gives the choice probability derived by the nested logit model in Equation (\ref{eqn:nl_2}). 
\end{theorem}

Theorem \ref{theorem:choice_equation_nl} proposes the use of entropy maximization formulation to derive probabilities that are equivalent to those obtained from the NL model. Proof can be found in Appendix \ref{appendix:nl}. To our knowledge, this is the first explicit incorporation of the NL model within the entropy maximization framework. A significant benefit of using the entropy maximization formulation instead of the NL model is the ability to set $\tau_B(m) = 0$. It is important to note that NL models operate within the domain of $\tau_B(m) > 0$ and cannot handle situations where $\tau_B(m) \rightarrow 0$. In contrast, the entropy maximization formulation provides a solution in cases where alternatives are perfectly correlated or purely relevant to each other, as exemplified by the well-known red bus and blue bus example.

\subsection{Utility maximization for hierarchical multinomial logit model}

\subsubsection{Hierarchical multinomial logit model}

We consider a simple hierarchical multinomial logit model that captures the destination and mode choice steps. Random utility and the Gumbel distribution are as follows. 
\begin{subequations}
\begin{align}
    &U_{ijm} = V_{ij} + V_{ijm} + E_{ij} + E_{ijm}  \\
    &F(\varepsilonvec_{ij}; \theta_j) = e^{- e^{- \theta_j \varepsilon_{ij}}}\\
    &F(\varepsilonvec_{ijm}; \theta_m) = e^{- e^{- \theta_m \varepsilon_{ijm}}}
\end{align}
\end{subequations}

\noindent Probability of choosing destination $j$ and mode $m$ can be expressed as follows. 

\begin{subequations}
\begin{align}
    & p_{j/i} 
    = \frac{\exp(\theta_j [V_{ij} + S_{ij}])}{\sum_{j'} \exp(\theta_{j} [V_{ij'} + S_{ij'}])}
    \label{eqn:prob_1}\\ 
    &p_{m/ij} 
    = \frac{\exp(\theta_m V_{ijm})}{\sum_{m'} \exp(\theta_{m} V_{ijm'})}
    \label{eqn:prob_2}\\
    & \text{where } S_{ij} = \frac{1}{\theta_m}\ln \sum_{m} \exp(\theta_m V_{ijm})
\end{align}
\end{subequations}

\subsubsection{Hierarchical multinomial logit model in entropy maximization framework}

\Problem{{\fontfamily{cmr}\selectfont HierMNL}}{$\overline{O_i}, \overline{T_{ij}}, \overline{T_{ijm}}, X^k_{ij}, X^q_{ijm}$}

\begin{subequations}
\begin{align}
    \max_{[p_{j/i}, p_{m/ij}]_{>0}} & - \sum_{ijm} \overline{O_i} P_{j/i} p_{m/ij} \ln (p_{j/i} p_{m/ij})\label{HeirMNL:objective} \\
    \text{s.t.} 
    & - \sum_{ij} \overline{O_i} p_{j/i} \ln (p_{j/j}) = 
    -\sum_{ij} \overline{T_{ij}} \ln(\frac{\overline{T_{ij}}}{\overline{O_{i}}}) && \quad [\frac{1}{\theta_j}]  \label{HeirMNL:ratio1} \\
    & - \sum_{ijm} \overline{O_i} P_{j/i} p_{m/ij} \ln (p_{m/ij}) = 
    -\sum_{ijm} \overline{T_{ijm}} \ln(\frac{\overline{T_{ijm}}}{\overline{T_{ij}}}) && \quad [\frac{1}{\theta_m}] \label{HeirMNL:ratio2} \\
    & \sum_{ij} \overline{O_i} p_{j/i} X^k_{ij} = \sum_{ij} \overline{T_{ij}} X^k_{ij}& , \forall k & \quad [\beta_k] \label{HeirMNL:aggregate1} \\
    & \sum_{ijm} \overline{O_i} p_{j/i} p_{m/ij} X^q_{ijm} = \sum_{ijm} \overline{T_{ijm}} X^q_{ijm}& , \forall q & \quad [\beta_q] \label{HeirMNL:aggregate2}\\
    & \sum_{j} p_{j/i} = 1 &,\forall i & \quad [\lambda_i] \label{HeirMNL:prob1}\\ 
    & \sum_{m} p_{m/ij} = 1 &,\forall i,j & \quad [\mu_{ij}]\label{HeirMNL:prob2}
\end{align}
\end{subequations}

Objective function in (\ref{HeirMNL:objective}) is the entropy function of $p_{m/i} \ln (p_{m/i})$, and $p_{m/i}$ is substituted with $p_{j/i} p_{m/ij}$. Constraint (\ref{HeirMNL:ratio1}) and (\ref{HeirMNL:ratio2}) ensure that the choice probability of each level can reproduce the ratio of observed trips. Constraint (\ref{HeirMNL:aggregate1}) and (\ref{HeirMNL:aggregate2}) state that the expectation of the aggregate value of each attribute should equal the observed aggregate
value. Constraint (\ref{HeirMNL:prob1}) and (\ref{HeirMNL:prob2}) ensure that the choice probabilities of a choice set should sum to unity. 

\begin{theorem}\label{theorem:hierMNL}
    In problem {\fontfamily{cmr}\selectfont HierMNL}, the optimal primal solutions $\widehat{p_{j/i}}, \widehat{p_{m/ij}}$ and the optimal dual solutions $\frac{1}{\widehat{\theta_m}}$, $\frac{1}{\widehat{\theta_j}}$, $\widehat{\beta_k}$, $\widehat{\beta_q}$ satisfy the folllowing hierarchical multinomial logit model. 
\end{theorem}

\begin{subequations}
\begin{align}
    &\widehat{p_{j/i} }
    = \frac{\exp(\widehat{\theta_j} [\sum_{k} \widehat{\beta_k} X^k_{ij} + S_{ij}])}{\sum_{j'} \exp(\widehat{\theta_{j}} [\sum_{k}\widehat{ \beta_k} X^k_{ij'} + S_{ij'}])}
    \label{hierMNL:prob_1}\\ 
    &\widehat{p_{m/ij}}
    = \frac{\exp(\widehat{\theta_m} [\sum_{q} \widehat{\beta_q} X^q_{ijm}])}{\sum_{m'} \exp(\widehat{\theta_{m}} [\sum_{q} \widehat{\beta_q} X^q_{ijm'}])}
    \label{hierMNL:prob_2}\\
    & \text{where } S_{ij} = \frac{1}{\widehat{\theta_m}}\ln \sum_{m} \exp(\widehat{\theta_m} \sum_{q} \widehat{\beta_q} X^q_{ijm})
\end{align}
\end{subequations}

\noindent where $S_{ij}$ is the satisfaction function, which is defined as the expected utility from a set of all modes. 

Theorem \ref{theorem:hierMNL} shows that the optimal solution of the entropy maximization formulation {\fontfamily{cmr}\selectfont HierMNL} is the same as the solution from the hierarchical multinomial logit model. This finding was first suggested by \citet{oppenheim1995urban}. We provide the full proof in Appendix \ref{appendix:hierMNL} for completeness. 

From Equation (\ref{hierMNL:prob_1}) and (\ref{hierMNL:prob_2}), it is evident that we can simply read off the dual variables $\beta_k, \beta_q$ from the optimal solution without a separate estimation process. \citet{anas1983discrete} demonstrates that the maximum likelihood estimation and entropy maximizing methods yield identical estimations for the MNL model by establishing the dual relationship between these two approaches. In addition, \citet{donoso2010microeconomic} provides a microeconomic interpretation of the entropy maximization problem, suggesting that the maximum entropy estimators of the multinomial logit model parameters reproduce rational user behavior. Moreover, \citet{donoso2011maximum} highlighted that the estimators from the entropy maximization are the better ones than the maximum likelihood estimators in the hierarchical logit model because they can replicate the observed fixed utility and the observed modal shares by incorporating additional constraints.

An interesting variant is when the scaling factors $\frac{1}{\theta_m}$, $\frac{1}{\theta_j}$, and parameters $\beta_k, \beta_q$ are given. The deterministic utilities are $\overline{V_{ij}} = \sum_{k} \overline{\beta_k} X^k_{ij}$, $\overline{V_{ijm}} = \sum_{q} \overline{\beta_q} X^q_{ijm}$. We can reconstruct an optimization problem from the Lagrangian function.

\Problem{{\fontfamily{cmr}\selectfont HierMNLVariant}}{$\overline{O_i}, \overline{T_{ij}}, \overline{T_{ijm}}, \overline{V_{ij}}, \overline{V_{ijm}}$, $\frac{1}{\theta_m}$, $\frac{1}{\theta_j}$}

\begin{subequations}
\begin{align}
    \max_{[p_{j/i}, p_{m/ij}]_{>0}} & - \frac{1}{\theta_j}[\sum_{ij} \overline{O_i}p_{j/i} ln(p_{j/i})]  - \frac{1}{\theta_m} [\sum_{ijm} \overline{O_i} p_{j/i} p_{m/ij} \ln(p_{m/ij})] \nonumber \\ & + \sum_{ij} \overline{O_i} p_{j/i}  \overline{V_{ij}}+ \sum_{ijm} \overline{O_i} p_{j/i} p_{m/ij} \overline{V_{ijm}} \\
    \text{s.t.} 
    & \sum_{j} p_{j/i} = 1 \quad \quad,\forall i & \quad [\lambda_i] \\ 
    & \sum_{m} p_{m/ij} = 1 \quad \quad,\forall i,j & \quad [\mu_{ij}]
\end{align}
\end{subequations}

Alternatively, instead of the decision variables $p_{j/i}$ and $p_{m/ij}$, $T_{ij}$ and $T_{ijm}$ can be adopted as previously established forms documented in the existing literature \citep{oppenheim1995urban, yang2009sensitivity}. 

\Problem{{\fontfamily{cmr}\selectfont HierMNLVariant2}}{$\overline{O_i}, \overline{T_{ij}}, \overline{T_{ijm}}, \overline{V_{ij}}, \overline{V_{ijm}}$, $\frac{1}{\theta_m}$, $\frac{1}{\theta_j}$}

\begin{subequations}
\begin{align}
    \max_{[T_{ij}, T_{ijm}]_{>0}} & - \frac{1}{\theta_{j'}}[\sum_{ij} T_{ij} ln(T_{ij})] - \frac{1}{\theta_m} [\sum_{ijm} T_{ijm} \ln(T_{ijm})]  + \sum_{ij} T_{ij}  \overline{V_{ij}}+ \sum_{ijm} T_{ijm} \overline{V_{ijm}} \\
    & \text{where } \frac{1}{\theta_{j'}} = \frac{1}{\theta_{j}} - \frac{1}{\theta_{m}} \nonumber \\
    \text{s.t.} 
    & \sum_{j} T_{ij} = \overline{O_i} \quad \quad,\forall i & \quad [\lambda_i]\\ 
    & \sum_{m} T_{ijm} = \overline{T_{ij}} \quad \quad,\forall i,j & \quad [\mu_{ij}] 
\end{align}
\end{subequations}

Lastly, it is important to note that the nested logit model and the hierarchical logit model are not the same, which can be confusing because they both incorporate hierarchical structures in modeling choices. The nested logit model focuses on the correlation among alternatives within nests, while the hierarchical logit model focuses on capturing choices made at different levels/stages in a hierarchical manner. To distinguish them more properly, we use the term hierarchical \textit{multinomial} logit model. The difference between the nested logit model and the hierarchical multinomial logit model is summarized in Table \ref{tab:nest_hier}.

\begin{table}[ht]
\caption{Comparison Between Nested Logit and Hierarchical Multinomial Logit Model}
\label{tab:nest_hier}
\begin{tabular}{|l|l|l|}
\hline
                              & \textbf{Nested logit}                                    & \textbf{Hierarchical multinomial logit }                  \\ \hline
\textbf{Random utility}       & \begin{tabular}[c]{@{}l@{}} $U_{m} = V_{m} + E_{m}$ \\ $F(\varepsilonvec; \tauvec) = e^{-\sum_{M \in \mathcal{M}} [\sum_{m \in M} e^{-\theta \varepsilon_m/{\tau_{M}}}]^{\tau_{M}}} $ \end{tabular} & \begin{tabular}[c]{@{}l@{}}$U_{ijm} = V_{ij} + V_{ijm} + E_{ij} + E_{ijm} $\\ $F(\varepsilonvec; \theta) = e^{-  e^{-\theta \varepsilon}}$\end{tabular} \\ \hline
\textbf{Probability}          & \begin{tabular}[c]{@{}l@{}}$p_M = \frac{e^{IV_{M}}}{\sum_{M' \in \mathcal{M}}e^{IV_{M'}}}$\\ $p_{m/M} =\frac{e^{\theta V_{m}/{\tau_{M}}}}{\sum_{m' \in M}e^{\theta V_{m'}/\tau_{M}}}$\\ $\text{where }IV_{\nest} = \tau_M \ln \sum_{m' \in \nest }e^{\theta V_{m'}/\tau_{\nest}}$ \end{tabular}  & \begin{tabular}[c]{@{}l@{}}$p_{j/i} 
    = \frac{\exp(\theta_j [V_{ij} + S_{ij}])}{\sum_{j'} \exp(\theta_{j} [V_{ij'} + S_{ij'}])}$\\ $p_{m/ij} 
    = \frac{\exp(\theta_m V_{ijm})}{\sum_{m'} \exp(\theta_{m} V_{ijm'})}$\\ $\text{where } S_{ij} = \frac{1}{\theta_m}\ln \sum_{m} \exp(\theta_m V_{ijm})$ \end{tabular}  \\ \hline
\end{tabular}
\end{table}

There is a similarity between inclusive value, $IV_M$, and satisfaction function, $S_{ij}$, as both of them represent aggregated values. $IV_M$ represents the attractiveness of nest $M$ as it measures the summation of normalized utilities of all alternatives in nest $M$. Similarly, $S_{ij}$ represents the expected utility of $ij$–pair from all choice alternatives in subsequent levels. 

\subsection{Hierarchical extended logit model for the combined travel demand modeling}

\subsubsection{Notation}
For the convenience of readers, we first present Table \ref{table:notations} which includes all the parameters and variables used in this document.

\begin{table}[ht]
\caption{Main Parameters and Decision Variables}\label{table:notations}
\begin{center}
\begin{tabular}{ll}
    \hline 
    $I$ & set of all origins \\
    \hline 
    $J$ & set of all destinations \\
    \hline 
    $\mathcal{M}$ & set of all modes \\
    \hline 
    $R$ & set of all routes \\
    \hline 
    $\theta_j, \theta_m, \theta_r$ & the scale factor for the Gumbel probability distribution of \\ &  the random utility error associated with destination, mode, and route\\
    \hline
    $\tau_{M}$ & dissimilarity factor of nest $M$ \\
    \hline
    $B(m)$ & function that identifies the nest to which alternative $m$ belongs \\
    \hline 
    $X^k_{ij}$ & the $k$-th attribute of $ij$-trip \\
    \hline
    $X^q_{ijm}$ & the $q$-th attribute of $ijm$-trip \\
    \hline
    $\beta_k, \overline{\beta_k}$ & parameter for $k$-th attribute \\
    \hline
    $\beta_q, \overline{\beta_q}$ & parameter for $q$-th attribute \\
    \hline
    $V_{ij}$ & the fixed utility for a traveler making a trip from origin $i$ to destination $j$ \\
    \hline
    $V_{ijm}$ & the fixed utility for a traveler making a trip from origin $i$ to destination $j$ \\& choosing mode $m$ \\
    \hline
    $V_{ijmr}$ & the fixed utility for a traveler making a trip from origin $i$ to destination $j$ \\& choosing mode $m$ and route $r$ \\
    \hline
    $g_{ijmr}$ & generalized travel cost of taking route $r$ on mode $m$  \\ & from origin $i$ to destination $j$ : $g_{ijmr} = \sum_{a_m \ni r} g^m_a(\cdot) \delta^{a_m}_{ijr}$ \\
    \hline
    $S_{ij}$ & the expected utility from a set of all mode $m \in \mathcal{M}$\\
    \hline 
    $S_{ijm}$ & the expected utility from a set of all route $r \in R$ \\
    \hline
    $IV_M$ & the inclusive value that represents the attractiveness of nest $M$, \\& or summation of normalized utilities of all alternatives in nest $M$ \\
    \hline
    $PS^{ijm}_r$ & path-size scaling parameter\\
    \hline 
    $\delta^{a_m}_{ijr}$ & link-route incidence indicator, 1 if link $a$ on the route $r$ \\& from origin $i$ to destination $j$ on mode $m$, 0 otherwise \\
    \hline
    $f^m_a, \overline{f^m_a}$ & flow in link $a$ using mode $m$ \\
    \hline
    $g^m_a(\cdot)$ & generalized link travel cost function for link $a$ of mode $m$ \\& : $g^m_a (\cdot) = \frac{1}{p} t^m_a(\cdot) + c^m_a$ \\&  where $\frac{1}{p}$ is value of time; $t^m_a(\cdot)$ is travel time function; $c^m_a$ travel cost\\
    \hline
    $g^m_{a}(\cdot)$ & function for mode specific link travel cost\\
    \hline
    $p_{j/i}$ & probability of choosing destination $j$ with the given origin $i$ \\
    \hline
    $p_{m/ij}$ & probability of choosing mode $m$ with given the origin $i$ and destination $j$\\
    \hline
    $p_{M/ij}$ & probability of choosing mode bundle $M$ with the given origin $i$ and destination $j$\\
    \hline
    $p_{m/M}$ & probability of choosing mode $m$ with the given mode bundle $M$\\
    \hline
    $p_{r/ijm}$ & probability of choosing route $r$ with the given origin $i$, destination $j$, and mode $m$\\
    \hline
    $\overline{O_i}$ & Given total demand from origin $i$\\
    \hline
    $T_{ij}$, $\overline{T_{ij}}$ & number of travelers from origin i to destination j, or a cell of a trip matrix \\
    \hline
    $T_{ijm}$, $\overline{T_{ijm}}$ & number of travelers using mode m from origin i to destination j \\
    \hline
    $T_{ijmr}$ & number of travelers taking route r on mode m from origin i to destination j \\
    \hline
\end{tabular}
\end{center}
{\raggedright \footnotesize We use the notation $\overline{\cdot}$ to denote that it is given value.\par}
\end{table}

\subsubsection{Hierarchical extended logit model}

\begin{figure}[H]
    \centering
    \includegraphics[width=13cm]{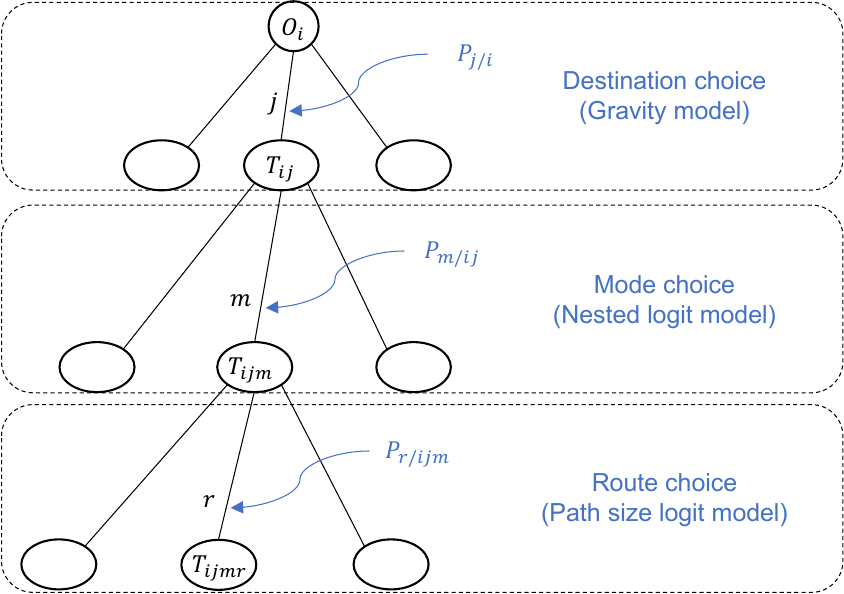}
    \caption{Hierarchical decision making process}
    \label{fig:hierarchical}
\end{figure}

Figure \ref{fig:hierarchical} represents a hierarchical structure of the individual traveler's decision making process. The decision tree is structured sequentially for the analytical purpose, but simultaneously solved in a combined model. This reflects choice decisions in reality as travelers consider all levels of choice simultaneously and these choices are interrelated with each other. The aggregation of rational individuals' destination, mode, and route choices correspond to collective utility maximization of the trip generation, trip distribution, modal split, and route assignment steps. 

Each step in the hierarchical tree is represented using an appropriate logit model. Destination choice, typically modeled with the gravity model, is represented by the multinomial logit model. As we discussed in section \ref{section:preliminary}, the entropy maximization modeling for the gravity model corresponds to that of the multinomial logit model. Mode choice is modeled with the nested logit model, which is often used to consider correlated alternatives in transportation modes. Route choice is modeled with the path-size logit model, which is developed by \citet{ben1999discrete} to handle route overlapping issues.

For the hierarchical structure, all levels of choices have to be considered simultaneously. When $\overline{O_i}$ is given, the choice probability $p_{j/i}, p_{m/ij}, p_{r/ijm}$ have to be jointly estimated. When the parameters associated with destination, mode, and route choice  $\theta_{j}$, $\theta_{m}$, $\theta_{r}$ are given, the total utility received from a single trip from origin $i$ to destination $j$ on mode $m$ and route $r$ can be specified as

\begin{subequations}
\begin{align}
    &U_{ijmr} = V_{ij} + V_{ijm} + V_{ijmr} + E_{ij} + E_{ijm} + E_{ijmr} \\
    &F(\varepsilonvec_{ij}; \theta_j) = e^{- e^{- \theta_j \varepsilon_{ij}}}\\
    &F(\varepsilonvec_{ijm}; \theta_m) = e^{- e^{- \theta_m \varepsilon_{ijm}}}\\
    &F(\varepsilonvec_{ijmr}; \theta_r) = e^{- e^{- \theta_r \varepsilon_{ijmr}}}
\end{align}
\end{subequations}

\noindent where $V_{ij}$, $V_{ijm}$, $V_{ijmr}$ are the fixed utilities for a traveler choosing destination $j$, mode $m$, and route $r$, and $\varepsilon_{ij}$, $\varepsilon_{ijm}$, $\varepsilon_{ijmr}$ are the corresponding error terms associated with each travel choice level. While $\varepsilon_{ij}$ and $\varepsilon_{ijmr}$ are defined as the Gumbel distribution as in Equation (\ref{eqn:gumbel_mnl}), the cumulative distribution function of error $\varepsilon_{ijm}$ embeds the correlation within nests as in Equation (\ref{eqn:NL_cdf}). 

The fixed utilities are defined as 

\begin{subequations}
    \begin{align}
        & V_{ij} = \sum_{k} \beta_k X^k_{ij}\\
        & V_{ijm} = \sum_{q} \beta_q X^q_{ijm}\\
        & V_{ijmr} = -g_{ijmr}
    \end{align}
\end{subequations}

\noindent where $X^k_{ij}$ and $X^q_{ijm}$ are observed travel attributes and $\beta_k$ and $\beta_q$ are parameters associated with the attributes. $g_{ijmr}$ is the generalized disutility (or travel cost) of taking route $r$ on mode $m$ from origin $i$ to destination $j$. 

The probability of choosing destination $p_{j/i}$, mode $p_{m/ij}$, and route $p_{r/ijm}$ can be written as follows.

\begin{subequations}
\begin{align}
    &p_{j/i} = 
     \frac{e^{ \theta_{j} ( V_{ij} + S_{ij})}}{\sum_{j' \in J} e^{\theta_{j} (V_{ij'}+S_{ij'}) }}\label{HierExtend:prob_1} \\
    &p_{m/ij} = p_{M/ij} p_{m/M}
    = \frac{e^{IV_{B(m)}}}{\sum_{M} e^{IV_{M}}} 
    \frac{e^{ \theta_{m} (V_{ijm} + S_{ijm} )/\tau_{B(m)}}}{\sum_{m' \in B(m)} e^{\theta_{m} (V_{ijm'} + S_{ijm'})/\tau_{B(m)}}}  \label{HierExtend:prob_2} \\
    & p_{r/ijm} 
    = \frac{PS^{ijm}_r e^{\theta_{r} V_{ijmr}}}{\sum_{r' \in B(r)} PS^{ijm}_{r'} e^{\theta_{r} V_{ijmr'}}} \label{HierExtend:prob_3} \\
    & \text{where }  S_{ij} = \frac{1}{\theta_m}\ln \sum_{M} e^{ IV_M} \label{HierExtend:S_ij} \\
    & \quad \quad IV_M = \tau_M \ln \sum_{m} e^{\theta_m [V_{ijm} + S_{ijm}] / \tau_M} \label{HierExtend:IV_M} \\
    & \quad \quad S_{ijm} = \frac{1}{\theta_r}\ln \sum_{r} PS^{ijm}_r e^{\theta_r V_{ijmr}} \label{HierExtend:S_ijm} \\
    & \quad \quad PS^{ijm}_r = \sum_{a_m \ni r} \frac{l_{a_m}}{L^{ijm}_r} \left( \frac{1}{\sum_{r'} \delta^{a_m}_{ijr'}} \right) \label{eqn:PS}
\end{align}
\end{subequations}

Equation (\ref{HierExtend:prob_1}) shows that the probability of choosing destination $j$ is affected by the fixed utility of choosing a destination, $V_{ij}$, and the expected utility from a set of all choices in the subsequent levels, $S_{ij}$. Equation (\ref{HierExtend:prob_2}) shows that the probability of choosing mode $m$ is expressed as a joint probability of choosing nest $M$, $p_{M/ij}$, and the probability of choosing alternative $m$ in the nest $M$, $p_{m/M}$. $p_{M/ij}$ is expressed by the inclusive value that represents the attractiveness of nest $M$, $IV_{M}$. $p_{m/M}$ is affected by the fixed utility of choosing a mode, $V_{ijm}$, and the expected utility from a set of all choices in the subsequent level, $S_{ijm}$. Equation (\ref{HierExtend:prob_3}) shows that the probability of choosing route $r$, $p_{r/ijm}$. Since route choice is the last step, the probability is affected only by the fixed utility of choosing the route, $V_{ijmr}$. Unlike the multinomial logit model, the path-size scaling factor, $PS^{ijm}_r$, is multiplied by the exponential utility term. Equation (\ref{HierExtend:S_ij}) and Equation (\ref{HierExtend:S_ijm}) show the satisfaction and Equation (\ref{HierExtend:IV_M}) shows the inclusive value, as we discussed in the previous section. As shown in Equation (\ref{eqn:PS}), each link $a_m$ in route $r$ is penalized according to the number of routes that share the link $a_m$, $\sum_{r'} \delta^{a_m}_{ijr'}$, and the significance of penalization is weighted according to how dominant link $a$ is in route $r$, $\sum_{a_m \ni r} \frac{l_{a_m}}{L^{ijm}_r} $. $\delta^{a_m}_{ijr}$ denotes the link-route incidence indicator, $l_{a_m}$ denotes the length of link $a_m$, and $L^{ijm}_r$ denotes length of route $r$.

Ultimately, given the positive number of travelers $\overline{O_i}$ at origin $i$, the number of travelers can be computed by multiplying the conditional probability at each stage in the hierarchical structure.

\begin{align}
    & T_{ij} = \overline{O_i} p_{j/i} \label{eqn:prob_destin}\\
    & T_{ijm} = T_{ij} p_{m/ij} \label{eqn:prob_mode}\\
    & T_{ijmr} = T_{ijm} p_{r/ijm} \label{eqn:prob_route}
\end{align}

Thus, the optimal flow pattern $\mathbb{T} = [T_{ij}, T_{ijm}, T_{ijmr}]$ from (\ref{eqn:prob_destin}) - (\ref{eqn:prob_route}) can be viewed as the equilibrium state because no traveler can improve one's utility by unilaterally changing the decision of destination/mode/route choice in the utility maximization framework.

\subsubsection{First stage – estimate parameters with observed data}

\Problem{{\fontfamily{cmr}\selectfont FirstStage}}{$\overline{O_i}$, $\overline{T_{ij}}$, $\overline{T_{ijm}}$, $\overline{T_{ijM}}$, $X^k_{ij}$, $X^q_{ijm}$, $PS^{ijm}_r$, $\frac{1}{\theta_r}$}

\begin{align}
    \max_{[p_{j/i}, p_{M/ij}, p_{m/M}, p_{r/ijm}]_{>0}} & 
    - \frac{1}{\theta_r}[\sum_{ijMm} \overline{O_i} p_{j/i} p_{M/ij} p_{m/M} p_{r/ijm} \ln (\frac{p_{r/ijm}}{PS^{ijm}_r})] \nonumber \\ 
    & - \sum_{m} \sum_{a_m} \int^{f^m_a}_{0} g^m_a(w)dw 
    \label{FirstStage:objective} 
\end{align}

\begin{subequations}
\begin{align}
    \text{s.t.} \nonumber 
    \\ & - \sum_{ij} \overline{O_i} p_{j/i} \ln (p_{j/j}) = 
    -\sum_{ij} \overline{T_{ij}} \ln(\frac{\overline{T_{ij}}}{\overline{O_{i}}}) && \quad [\frac{1}{\theta_j}]  \label{FirstStage:ratio1} \\
    & - \sum_{ijM} \overline{O_i} p_{j/i} p_{M/ij} \ln (p_{M/ij}) = 
    -\sum_{ijm} \overline{T_{ijM}} \ln(\frac{\overline{T_{ijM}}}{\overline{T_{ij}}}) && \quad [\frac{1}{\theta_m}] \label{FirstStage:ratio2}\\
    & - \sum_{ijm} \overline{O_i} p_{j/i} p_{M/ij} p_{m/M} \ln (p_{m/M}) = 
     - \sum_{ijm} \overline{T_{ijm}} \ln(\frac{\overline{T_{ijm}}}{\overline{T_{ijM}}}) &, \forall M & \quad [\frac{\tau_M}{\theta_m}] \label{FirstStage:ratio3}\\
    & \sum_{ij} \overline{O_i} p_{j/i} X^k_{ij} = \sum_{ij} \overline{T_{ij}} X^k_{ij}& , \forall k & \quad [\beta_k] \label{FirstStage:aggregate1} \\
    & \sum_{ijm} \overline{O_i} p_{j/i} p_{M/ij} p_{m/M} X^q_{ijm} = \sum_{ijm} \overline{T_{ijm}} X^q_{ijm} & , \forall q & \quad [\beta_q] \label{FirstStage:aggregate2}\\
    & \sum_{j} p_{j/i} = 1 &,\forall i & \quad [\lambda_i]\label{FirstStage:prob1}\\ 
    & \sum_{M} p_{M/ij} = 1 &,\forall i,j & \quad [\mu_{ij}]\label{FirstStage:prob2}\\
    & \sum_{m} p_{m/M} = 1 &,\forall M & \quad [\kappa_{M}]\label{FirstStage:prob3} \\
    & \sum_{r} p_{r/ijm} = 1 &,\forall r & \quad [\nu_{ijm}]\label{FirstStage:prob4} \\
    & f^m_a = \sum_{ijr} \overline{O_i} p_{j/i} p_{M/ij} p_{m/M} p_{r/ijm} \delta^{a_m}_{ijr} \label{FirstStage:flow}
\end{align}
\end{subequations}

Lagrangian multipliers are in the square brackets. The first term of the objective function (\ref{FirstStage:objective}) is obtained from the entropy function $- p_{r/i} \ln (p_{r/i})$, substituting $(p_{r/i})$ by $p_{j/i} p_{M/ij} p_{m/M} p_{r/ijm}$ and removing the constant terms using Equation (\ref{FirstStage:ratio1})–(\ref{FirstStage:ratio3}). The second term is the form of the Beckmann equation contributing to finding user equilibrium in route choice. Constraint (\ref{FirstStage:ratio1}) – (\ref{FirstStage:ratio3}) restrict entropy in each group. Constraint (\ref{FirstStage:aggregate1}) and (\ref{FirstStage:aggregate2}) ensure that the fixed utility and the market shares are aligned with observed values. Constraint (\ref{FirstStage:prob1})–(\ref{FirstStage:prob4}) ensure that the choice probabilities within a choice set sum up to one. Equation (\ref{FirstStage:flow}) shows the relationship between link flow and route flow using the indicator $\delta^{a_m}_{ijr}$.  


\begin{assumption}\label{assume:costfunc}
    The generalized cost function $g^m_a(\mathbb{f})$ is continuous, differentiable, strictly increasing, and separable (i.e., they only rely on their own flows $g^m_a(\mathbb{f})=g^m_a(f^m_a)$). Also, the generalized route cost function $g_{ijmr}$ is addictive, i.e., $g_{ijmr} = \sum_{a_{m}} g^m_a (\cdot) \delta^{a_m}_{ijr}$.
\end{assumption}

\begin{theorem}\label{theorem:HierExtend}
In problem {\fontfamily{cmr}\selectfont FirstStage}, the optimal primal solutions $\widehat{p_{j/i}}$, $\widehat{p_{M/ij}}, \widehat{p_{m/M}}$, $\widehat{p_{r/ijm}}$ and the optimal dual solution $\frac{1}{\widehat{\theta_j}}$, $\frac{1}{\widehat{\theta_m}}$, $\frac{\widehat{\tau_M}}{\widehat{\theta_m}}$, $\widehat{\beta_k}$, $\widehat{\beta_q}$ satisfy the hierarchical extended logit model defined in (\ref{HierExtend:prob_1}) - (\ref{eqn:PS}) 
\end{theorem}

Under Assumption \ref{assume:costfunc}, the equivalent minimization problem of {\fontfamily{cmr}\selectfont FirstStage} is a convex programming with a convex objective function and convex constraints. Furthermore, Theorem \ref{theorem:HierExtend} highlights that the utility maximization formulation in problem {\fontfamily{cmr}\selectfont FirstStage} presents an alternative approach to the hierarchical extended logit model. Proof of Theorem \ref{theorem:HierExtend} can be found in Appendix \ref{appendix:hier_extend}.

\Problem{{\fontfamily{cmr}\selectfont FirstStageVariant}}{$\overline{O_i}, \overline{T_{ij}}, \overline{T_{ijm}}, \overline{T_{ijM}}, X^k_{ij}, X^q_{ijm}$, $PS^{ijm}_r$, $\frac{1}{\theta_r}$, $\sigma, \overline{f^m_a}$}

\begin{align}
    & \max_{[p_{j/i}, p_{M/ij}, p_{m/M}, p_{r/ijm}, f^m_a]_{>0}} \text{Equation (\ref{FirstStage:objective})} + \sigma \sum_{m} \sum_{a_m} |f^m_a - \overline{f^m_a}| \\
    & \text{s.t. Constraint (\ref{FirstStage:ratio1})–(\ref{FirstStage:flow})}
\end{align}

The useful variant is created by adding the error between output link flow values and observed values. By adding the additional term, the optimal solution of problem {\fontfamily{cmr}\selectfont FirstStageVariant} is no longer the user equilibrium state. Since the real world does not exhibit a perfect user equilibrium, we need to balance between user equilibrium and observation. $\sigma$ is the balancing parameter.

\subsubsection{Second stage – estimate future travel with the fixed parameters}

Now, future estimation of $\widetilde{O_i}$ is given. Travel attributes in the new transportation network, $\widetilde{X^k_{ij}}$ and $ \widetilde{X^q_{ijm}}$, are also given. $\widehat{\theta_j}$, $\widehat{\theta_m}$, $\widehat{\tau_M}$, $\widehat{\beta_k}$, $\widehat{\beta_q}$ are estimated in the first stage. Accordingly, the deterministic utilities are $\overline{V_{ij}} = \sum_{k} \widehat{\beta_k} \widetilde{X^k_{ij}}$, $\overline{V_{ijm}} = \sum_{q} \widehat{\beta_q} \widetilde{X^q_{ijm}}$. We can reconstruct an optimization problem from the Lagrangian function.

\Problem{{\fontfamily{cmr}\selectfont SecondStage}}{$\widetilde{O_i}, \overline{T_{ij}}, \overline{T_{ijm}}, \overline{T_{ijM}}, \overline{V_{ij}}, \overline{V_{ijm}}$, $\frac{1}{\theta_r}$, $\frac{1}{\widehat{\theta_m}}$, $\frac{1}{\widehat{\theta_j}}$, $\widehat{\tau_M}$, $PS^{ijm}_r$}

\begin{subequations}
\begin{align}
    \max_{[p_{j/i}, p_{M/ij}, p_{m/M}, p_{r/ijm}]_{>0}} 
    & 
    - \frac{1}{\widehat{\theta_j}} [\sum_{ij} \widetilde{O_i} P_{j/i} \ln (p_{j/j})] \\
    & -\frac{1}{\widehat{\theta_m}} \left[ \sum_{ijMm} \widetilde{O_i} P_{j/i} p_{M/ij} p_{m/M} \ln (p_{M/ij} p_{m/M}) \right. \\
    & \quad \quad \left. + \widehat{\tau_{M}} [\sum_{ijm} \widetilde{O_i} P_{j/i} p_{M/ij} p_{m/M} \ln (p_{m/M}) ] \right] \\
    & - \frac{1}{\theta_r}[\sum_{ijm} \widetilde{O_i} P_{j/i} p_{M/ij} p_{m/M} p_{r/ijm} \ln (\frac{p_{r/ijm}}{PS^{ijm}_r})] \nonumber \\ 
    & - \sum_{m} \sum_{a_m} \int^{f^m_a}_{0} g^m_a(w)dw  \\
    & + \sum_{ij} \widetilde{O_i} p_{j/i} \overline{V_{ij}} + \sum_{ijm} \widetilde{O_i} p_{j/i} p_{m/ij} \overline{V_{ijm}}
\end{align}
\end{subequations}

\begin{subequations}
\begin{align}
    \text{s.t.} \\ 
    & \sum_{j} p_{j/i} = 1 &,\forall i & \quad [\lambda_i]\label{SecondStage:prob1}\\ 
    & \sum_{M} p_{M/ij} = 1 &,\forall i,j & \quad [\mu_{ij}]\label{SecondStage:prob2}\\
    & \sum_{m} p_{m/M} = 1 &,\forall M & \quad [\kappa_{M}]\label{SecondStage:prob3} \\
    & \sum_{r} p_{r/ijm} = 1 &,\forall r & \quad [\nu_{r}]\label{SecondStage:prob4} \\
    & f^m_a = \sum_{ijr} \widetilde{O_i} p_{j/i} p_{M/ij} p_{m/M} p_{r/ijm} \delta^{a_m}_{ijr} \label{SecondStage:flow}
\end{align}
\end{subequations}

\begin{corollary}\label{theorem:SecondStage}
    In problem {\fontfamily{cmr}\selectfont SecondStage}, the optimal solutions $\widehat{p_{j/i}}$, $\widehat{p_{M/ij}}, \widehat{p_{m/M}}$, $\widehat{p_{r/ijm}}$ satisfy the hierarchical extended logit model defined in (\ref{HierExtend:prob_1}) - (\ref{eqn:PS}). 
\end{corollary}

Under Assumption \ref{assume:costfunc}, the equivalent minimization problem of {\fontfamily{cmr}\selectfont SecondStage} is convex programming with the convex objective function and linear constraints. Corollary \ref{theorem:SecondStage} supports that the future estimation can be obtained by solving problem {\fontfamily{cmr}\selectfont SecondStage}. From the optimal solution $\widehat{p_{j/i}}$, $\widehat{p_{M/ij}}, \widehat{p_{m/M}}$, $\widehat{p_{r/ijm}}$, we get the future estimated trip distribution $T_{ij}, T_{ijm}, T_{ijmr}, f^m_a$ from Equation (\ref{eqn:prob_destin})-(\ref{eqn:prob_route}) and (\ref{SecondStage:flow}).

In summary, the {\fontfamily{cmr}\selectfont FirstStage} incorporates observed data from each decision-making step as constraints and solves the problem simultaneously, resulting in consistent outcomes across all levels. Moreover, this formulation gains greater versatility as it allows the introduction of a balancing term between user equilibrium in the utility maximization framework and observed link flow values. This recognizes real-world conditions, which may not perfectly align with user equilibrium, leading to a more accurate representation of the actual situation. As a result, the {\fontfamily{cmr}\selectfont SecondStage} model is enhanced by incorporating fixed parameters that better reflect reality, leading to more reliable forecasts.

\section{Conclusion}

The travel demand forecasting model plays a critical role in assessing large-scale infrastructure projects. Sequential four-step travel demand modeling remains widely used despite the inherent limitations, i.e., the absence of a unifying framework and discrepancies between input and output values in each step. Although combined models have been developed, they have not been practically applicable due to a lack of satisfaction in both behavioral richness and computational efficiency, which are desirable in the real world. To bridge this gap, we propose a convex programming approach equivalent to the hierarchical extended logit model. The most significant advantage of our method lies in its ability to be modeled as convex programming, ensuring solution existence, uniqueness, and tractability.

Our model is built upon the utility maximization framework, leveraging well-established findings. Specifically, an individual choice of destination, mode, and route can be expressed as random utility, and the trip generation/distribution, modal split, and traffic assignment can be viewed as maximizing the aggregated utility. Building upon that, we extend the well-known equivalence of the multinomial logit model and the satisfaction maximization problem to the nested logit setting. This extension allows our model to exhibit greater behavioral richness compared to previous convex programming approaches. Each step is expressed using appropriate models, such as the multinomial logit model (which is equivalent to the gravity model) for destination choice, the nested logit model for mode choice, and the path-size logit model for route choice.

The combined model demonstrates its advantages during the two stages of data assimilation and forecasting. In the data assimilation stage, it serves as a viable alternative to maximum likelihood estimation by incorporating observations as constraints. This approach maximizes the log-likelihood of observing the current trip distribution while preserving the observed share of each choice. Unlike traditional models with feedback loops, our convex programming solution allows for the simultaneous estimation of all travel levels, thereby overcoming the input-output discrepancies in each level. Moreover, the combined model's versatility enables us to strike a balance between theoretical user equilibrium and real-world observations (e.g., traffic counts at each link). By achieving a better representation of reality, we obtain more reliable future estimations.

In future research, exploring the application of this model to various transportation networks in different cities would provide highly valuable insights. Additionally, conducting a comparative analysis between the performance of the convex programming model and traditional four-step sequential approaches, as well as other advanced travel demand forecasting techniques like agent-based modeling or machine learning approaches, could offer valuable assessments of its strengths and limitations.

For a more comprehensive understanding of travelers' behavior, there are promising extensions to consider. Firstly, accommodating multi-modal trip chains allows for a more comprehensive capture of travelers' behavior. Secondly, investigating dynamic travel demand modeling, where travel patterns evolve over time would better reflect real-world dynamics. To address the challenge of non-separable link cost functions, where the link travel time of each mode depends on the flows of other modes on that link, further research could explore effective ways to handle these complexities, thereby expanding the model's applicability.

The extended model's capability to capture travelers' responses to different policy interventions, such as tolling, holds significant potential for achieving sustainable urban mobility goals. Additionally, the model can be utilized to evaluate the investment in emerging technologies, such as autonomous vehicles and electric vehicles. By exploring these extensions and applications, the proposed convex programming model can evolve into a powerful and versatile tool, aiding transportation planning and policy-making in diverse urban contexts.


\bibliographystyle{unsrtnat}
\bibliography{references}  

\begin{thebibliography}{30}
\providecommand{\natexlab}[1]{#1}
\providecommand{\url}[1]{\texttt{#1}}
\expandafter\ifx\csname urlstyle\endcsname\relax
  \providecommand{\doi}[1]{doi: #1}\else
  \providecommand{\doi}{doi: \begingroup \urlstyle{rm}\Url}\fi

\bibitem[de~Dios~Ort{\'u}zar and Willumsen(2011)]{de2011modelling}
Juan de~Dios~Ort{\'u}zar and Luis~G Willumsen.
\newblock \emph{Modelling transport}.
\newblock John wiley \& sons, 2011.

\bibitem[McNally(2007)]{mcnally2007four}
Michael~G McNally.
\newblock The four-step model.
\newblock In \emph{Handbook of transport modelling}, volume~1, pages 35--53.
  Emerald Group Publishing Limited, 2007.

\bibitem[Boyce(2002)]{boyce2002sequential}
David Boyce.
\newblock Is the sequential travel forecasting paradigm counterproductive?
\newblock \emph{Journal of urban planning and development}, 128\penalty0
  (4):\penalty0 169--183, 2002.

\bibitem[Garrett and Wachs(1996)]{garrett1996transportation}
Mark Garrett and Martin Wachs.
\newblock \emph{Transportation planning on trial: The Clean Air Act and travel
  forecasting}.
\newblock Sage Publications, 1996.

\bibitem[Oppenheim et~al.(1995)]{oppenheim1995urban}
Norbert Oppenheim et~al.
\newblock \emph{Urban travel demand modeling: from individual choices to
  general equilibrium.}
\newblock John Wiley and Sons, 1995.

\bibitem[Yao et~al.(2014)Yao, Chen, Ryu, and Shi]{yao2014general}
Jia Yao, Anthony Chen, Seungkyu Ryu, and Feng Shi.
\newblock A general unconstrained optimization formulation for the combined
  distribution and assignment problem.
\newblock \emph{Transportation Research Part B: Methodological}, 59:\penalty0
  137--160, 2014.

\bibitem[Zhou et~al.(2009)Zhou, Chen, and Wong]{zhou2009}
Zhong Zhou, Anthony Chen, and Sze~Chun Wong.
\newblock Alternative formulations of a combined trip generation, trip
  distribution, modal split, and trip assignment model.
\newblock \emph{European Journal of Operational Research}, 198\penalty0
  (1):\penalty0 129--138, 2009.

\bibitem[Shannon(1948)]{shannon1948mathematical}
Claude~Elwood Shannon.
\newblock A mathematical theory of communication.
\newblock \emph{The Bell system technical journal}, 27\penalty0 (3):\penalty0
  379--423, 1948.

\bibitem[Fang et~al.(1997)Fang, Rajasekera, and Tsao]{fang1997entropy}
Shu-Cherng Fang, Jay~R Rajasekera, and H-S~Jacob Tsao.
\newblock \emph{Entropy optimization and mathematical programming}, volume~8.
\newblock Springer Science \& Business Media, 1997.

\bibitem[Wilson(1969)]{wilson1969use}
Alan~Geoffrey Wilson.
\newblock The use of entropy maximising models, in the theory of trip
  distribution, mode split and route split.
\newblock \emph{Journal of transport economics and policy}, pages 108--126,
  1969.

\bibitem[McFadden(1973)]{mcfadden1973conditional}
Daniel McFadden.
\newblock Conditional logit analysis of qualitative choice behaviour.
\newblock In P.~Zarembka, editor, \emph{Frontiers in Econometrics}. Academic
  Press, New York, 1973.

\bibitem[Ben-Akiva et~al.(1985)Ben-Akiva, Lerman, Lerman,
  et~al.]{ben1985discrete}
Moshe~E Ben-Akiva, Steven~R Lerman, Steven~R Lerman, et~al.
\newblock \emph{Discrete choice analysis: theory and application to travel
  demand}, volume~9.
\newblock MIT press, 1985.

\bibitem[Miyagi and Morisugi(1996)]{miyagi1996direct}
Toshihiko Miyagi and Hisa Morisugi.
\newblock A direct measure of the value of choice-freedom.
\newblock \emph{Papers in Regional Science}, 75\penalty0 (2):\penalty0
  121--134, 1996.

\bibitem[Anas(1983)]{anas1983discrete}
Alex Anas.
\newblock Discrete choice theory, information theory and the multinomial logit
  and gravity models.
\newblock \emph{Transportation Research Part B: Methodological}, 17\penalty0
  (1):\penalty0 13--23, 1983.

\bibitem[Beckmann et~al.(1956)Beckmann, McGuire, and
  Winsten]{beckmann1956studies}
Martin Beckmann, Charles~B McGuire, and Christopher~B Winsten.
\newblock \emph{Studies in the Economics of Transportation}.
\newblock YALE UNIVERSITY PRESS, 1956.

\bibitem[Evans(1976)]{evans1976derivation}
Suzanne~P Evans.
\newblock Derivation and analysis of some models for combining trip
  distribution and assignment.
\newblock \emph{Transportation research}, 10\penalty0 (1):\penalty0 37--57,
  1976.

\bibitem[Florian et~al.(1975)Florian, Nguyen, and Ferland]{florian1975combined}
Michael Florian, Sang Nguyen, and Jacques Ferland.
\newblock On the combined distribution-assignment of traffic.
\newblock \emph{Transportation Science}, 9\penalty0 (1):\penalty0 43--53, 1975.

\bibitem[Oppenheim(1993)]{oppenheim1993equilibrium}
Norbert Oppenheim.
\newblock Equilibrium trip distribution/assignment with variable destination
  costs.
\newblock \emph{Transportation Research Part B: Methodological}, 27\penalty0
  (3):\penalty0 207--217, 1993.

\bibitem[Florian(1977)]{florian1977traffic}
Michael Florian.
\newblock A traffic equilibrium model of travel by car and public transit
  modes.
\newblock \emph{Transportation Science}, 11\penalty0 (2):\penalty0 166--179,
  1977.

\bibitem[Abdulaal and LeBlanc(1979)]{abdulaal1979methods}
Mustafa Abdulaal and Larry~J LeBlanc.
\newblock Methods for combining modal split and equilibrium assignment models.
\newblock \emph{Transportation Science}, 13\penalty0 (4):\penalty0 292--314,
  1979.

\bibitem[Fern{\'a}ndez et~al.(1994)Fern{\'a}ndez, De~Cea, Florian, and
  Cabrera]{fernandez1994network}
Enrique Fern{\'a}ndez, Joaquin De~Cea, Michael Florian, and Enrique Cabrera.
\newblock Network equilibrium models with combined modes.
\newblock \emph{Transportation Science}, 28\penalty0 (3):\penalty0 182--192,
  1994.

\bibitem[Garc{\'i}a and Mar{\'i}n(2005)]{garcia2005network}
Ricardo Garc{\'i}a and Angel Mar{\'i}n.
\newblock Network equilibrium with combined modes: models and solution
  algorithms.
\newblock \emph{Transportation Research Part B: Methodological}, 39\penalty0
  (3):\penalty0 223--254, 2005.

\bibitem[Florian and Nguyen(1978)]{florian1978}
Michael Florian and Sang Nguyen.
\newblock A combined trip distribution modal split and trip assignment mode.
\newblock \emph{Transportation Research}, 12\penalty0 (4):\penalty0 241--246,
  1978.

\bibitem[Friesz(1981)]{friesz1981equivalent}
Terry~L Friesz.
\newblock An equivalent optimization problem for combined multiclass
  distribution, assignment and modal split which obviates symmetry
  restrictions.
\newblock \emph{Transportation Research Part B: Methodological}, 15\penalty0
  (5):\penalty0 361--369, 1981.

\bibitem[Safwat and Magnanti(1988)]{safwat1978}
K.~Nabil~Ali Safwat and Thomas~L. Magnanti.
\newblock A combined trip generation, trip distribution, modal split, and trip
  assignment model.
\newblock \emph{Transportation Science}, 22\penalty0 (1):\penalty0 14--30,
  1988.

\bibitem[Wen and Koppelman(2001)]{wen2001generalized}
Chieh-Hua Wen and Frank~S Koppelman.
\newblock The generalized nested logit model.
\newblock \emph{Transportation Research Part B: Methodological}, 35\penalty0
  (7):\penalty0 627--641, 2001.

\bibitem[Donoso and de~Grange(2010)]{donoso2010microeconomic}
Pedro Donoso and Louis de~Grange.
\newblock A microeconomic interpretation of the maximum entropy estimator of
  multinomial logit models and its equivalence to the maximum likelihood
  estimator.
\newblock \emph{Entropy}, 12\penalty0 (10):\penalty0 2077--2084, 2010.

\bibitem[Donoso et~al.(2011)Donoso, De~Grange, and
  Gonz{\'a}lez]{donoso2011maximum}
Pedro Donoso, Louis De~Grange, and Felipe Gonz{\'a}lez.
\newblock A maximum entropy estimator for the aggregate hierarchical logit
  model.
\newblock \emph{Entropy}, 13\penalty0 (8):\penalty0 1425--1445, 2011.

\bibitem[Yang and Chen(2009)]{yang2009sensitivity}
Chao Yang and Anthony Chen.
\newblock Sensitivity analysis of the combined travel demand model with
  applications.
\newblock \emph{European Journal of Operational Research}, 198\penalty0
  (3):\penalty0 909--921, 2009.

\bibitem[Ben-Akiva and Bierlaire(1999)]{ben1999discrete}
Moshe Ben-Akiva and Michel Bierlaire.
\newblock Discrete choice methods and their applications to short term travel
  decisions.
\newblock In \emph{Handbook of transportation science}, pages 5--33. Springer,
  1999.

\end{thebibliography}

\section{Appendix}

\subsection{Proof for entropy maximizing model for nested logit model}\label{appendix:nl}

\begin{proof}

Let $p_M$ and $p_{m/M}$ be independent variables where $p_m = p_M p_{m/M}$. From Equation (\ref{eqn:max_choice_equation_nl}), convert the problem to minimization problem, multiply $\theta$ in the objective function, and then substitute $p_m$ with $p_M p_{m/M}$. Then, we obtain the equivalent problem as below. 

\begin{subequations}\label{eqn:Z_function}
\begin{align}
    & \min_{p_M, p_{m/M}} Z 
    = \sum_{M \in \mathcal{M}} p_M \ln p_M 
    + \sum_{M \in \mathcal{M}} \tau_M p_M (\sum_{m \in M} \left(  p_{m/M} \ln p_{m/M} \right)) 
    - \sum_{M \in \mathcal{M}} \sum_{m \in M} p_{M} p_{m/M} \theta V_m
    \tag{\ref{eqn:Z_function}}
\end{align}
\begin{align}
    \text{subject to } & \sum_{M \in \mathcal{M}} p_M = 1  & [\lambda] \label{eqn:b_sum}\\
    & \sum_{m \in M} p_{m/M} = 1, \quad \forall M \in \mathcal{M} & [\varphi_{M}] \label{eqn:k_sum}
\end{align}
\end{subequations}

\noindent First order condition:
\begin{align}
    \frac{\partial Z}{\partial p_M} &= \ln p_M + 1 + \tau_M \sum_{m \in M} (p_{m/M}\ln(p_{m/M})) - \sum_{m \in M} p_{m/M} \theta V_m - \lambda = 0 \label{proofnl:partial_1}\\
    \frac{\partial Z}{\partial p_{m/M}} &= \tau_{M} p_M (\ln p_{m/M} + 1) - p_M \theta V_m - \varphi_M =0 \label{eqn:partial_2}
\end{align}

\noindent We first express $p_{m/M}$ as $p_M$ using Equation (\ref{eqn:partial_2}). 

\begin{align}
    &\ln p_{m/M} = \frac{p_M \theta V_m + \varphi_M}{\tau_M p_M} - 1 \nonumber\\
    & \Rightarrow p_{m/M} = \frac{\exp(\theta V_{m}/\tau_M)}{\exp(1-\varphi_M/\tau_M p_M)} \label{proofnl:p_mM} \\
    & \text{Using} \sum_{m \in M} p_{m/M} = 1 \text{ in Constraint (\ref{eqn:k_sum}),} \nonumber \\
    & \Rightarrow \sum_{m \in M} \frac{\exp(\theta V_{m}/\tau_M)}{\exp(1-\varphi_M/\tau_M p_M)} = 1  \nonumber \\
    &\text{The terms that are not related to $m$ cancels out.} \nonumber \\
    & \Rightarrow \sum_{m \in M} \exp(\theta V_{m}/\tau_M) = \exp(1-\varphi_M/\tau_M p_M) \label{proofnl:denominator} \\
    & \text{Using Equation (\ref{proofnl:denominator}), Equation (\ref{proofnl:p_mM}) can be rewritten. } \nonumber \\
    & \therefore p_{m/M} = \frac{\exp(\theta V_m/\tau_{M})}{\sum_{m' \in M} \exp(\theta V_{m'}/\tau_{M})} \label{proofnl:p_mM_final}
\end{align}

Using the definition, $IV_{M} = \tau_{M} \ln \sum_{m' \in M} \exp({\theta V_{m'}/\tau_{M}})$, Equation (\ref{proofnl:p_mM_final}) can be rewritten. 

\begin{align}
    &p_{m/M} = \exp(\frac{\theta V_m}{\tau_M} - \frac{IV_M}{\tau_M}) \nonumber \\
    & \Rightarrow \ln p_{m/M} = \frac{\theta V_m}{\tau_M} - \frac{IV_M}{\tau_M} \label{proofnl:ln_pkb}
\end{align}

\noindent Now express $p_M$ as $p_{m/M}$ using Equation (\ref{proofnl:partial_1}) and substitute $\ln p_{m/M}$ by Equation (\ref{proofnl:ln_pkb}).

\begin{align*}
    & \ln p_M + 1 + \tau_M \sum_{m \in M} p_{m/M} \left( \frac{\theta V_m}{\tau_M} - \frac{IV_M}{\tau_M} \right) - \sum_{m \in M} p_{m/M} \theta V_{m} - \lambda = 0 \\
    & \Rightarrow \ln p_M = IV_M + \lambda - 1 \\
    & \text{Similarly, using Constraint (\ref{eqn:b_sum}), we get} 
\end{align*}

\begin{align}
    & p_M = \frac{ \exp(IV_M)}{\sum_{M \in \mathcal{M}} \exp(IV_M)}
\end{align}

Thus, the optimal solution of the problem equivalent to the probability expressed as a closed form nested logit model. 

\begin{align}
    p_{M} p_{m/M} = \frac{ \exp(IV_M)}{\sum_{M' \in \mathcal{M}} \exp(IV_{M'})} \frac{\exp(\theta V_m/\tau_{M})}{\sum_{m' \in M} \exp(\theta V_{m'}/\tau_{M})}
\end{align}

\end{proof}

\subsection{Proof for hierarhical multinomial logit model}\label{appendix:hierMNL}

\begin{proof}
    Objective function in Equation (\ref{HeirMNL:objective}) can be rewritten as below.
\begin{subequations}
    \begin{align*}
        & - \sum_{ijm} \overline{O_i} p_{j/i} p_{m/ij} \ln (p_{j/i} p_{m/ij}) \\
        & = - \sum_{ij} \overline{O_i} p_{j/i} \ln (p_{j/i}) - \sum_{ijm} \overline{O_i} p_{j/i} p_{m/ij} \ln (p_{m/ij}) 
    \end{align*}
\end{subequations}

From Constraint (\ref{HeirMNL:ratio1}) and (\ref{HeirMNL:ratio2}), the objective function can be treated as a constant, thus we remove the constant term in the below Lagrangian function. 

\begin{align*}
    \mathcal{L} = & - \frac{1}{\theta_j} [\sum_{ij} \overline{O_i} p_{j/i} \ln (p_{j/i}) 
    -\sum_{ij} \overline{T_{ij}} \ln(\frac{\overline{T_{ij}}}{\overline{O_{i}}})] \\
    & - \frac{1}{\theta_m} [\sum_{ijm} \overline{O_i} p_{j/i} p_{m/ij} \ln (p_{m/ij})  
    -\sum_{ijm} \overline{T_{ijm}} \ln(\frac{\overline{T_{ijm}}}{\overline{T_{ij}}}) ] \\
    & + \sum_{k} \beta_k [\sum_{ij} \overline{O_i} p_{j/i} X^k_{ij} - \sum_{ij} \overline{T_{ij}} X^k_{ij}]  \\
    & + \sum_{q} \beta_q [\sum_{ijm} \overline{O_i} p_{j/i} p_{m/ij} X^q_{ijm} - \sum_{ijm} \overline{T_{ijm}} X^q_{ijm}] \\
    & + \sum_{i} \lambda_i [\sum_{j} p_{j/i} - 1] \\
    & + \sum_{ij} \mu_{ij} [\sum_{m} p_{m/ij} - 1] \\
\end{align*}

\noindent Write the first order conditions.
\begin{subequations}
    \begin{align}
        & \frac{\partial \mathcal{L}}{\partial p_{j/i}} = -\frac{1}{\theta_j} [\overline{O_i} (\ln p_{j/i} + 1)]   -\frac{1}{\theta_m} [\overline{O_i} \sum_{m}p_{m/ij} \ln p_{m/ij}]
        + \overline{O_i}  \sum_{k} \beta_k X^k_{ij} \nonumber \\
        & \quad \quad \quad + \overline{O_i} \sum_{m} p_{m/ij} \sum_{q} \beta_q X^q_{ijm} + \lambda_i = 0 \label{eqn:1}\\
        & \frac{\partial \mathcal{L}}{\partial p_{m/ij}} = -\frac{1}{\theta_m} [\overline{O_i} p_{j/i} (\ln p_{m/ij} + 1)] + \overline{O_i} p_{j/i} \sum_{q} \beta_q X^q_{ijm} + \mu_{ij}  = 0 \label{eqn:2}
    \end{align}
\end{subequations}

From Equation (\ref{eqn:2}), we obtain the equation below by dividing it into $\overline{O_i} p_{j/i}$. 

\begin{subequations}
    \begin{align}
   & \Rightarrow -\frac{1}{\theta_m} [ (\ln p_{m/ij} + 1)]  + \sum_{q} \beta_q X^q_{ijm} + \frac{\mu_{ij}}{\overline{O_i} p_{j/i}} = 0 \nonumber \\ 
   & \Rightarrow \ln p_{m/ij} = \theta_m \sum_{q} \beta_q X^q_{ijm} + \theta_m  \frac{\mu_{ij}}{\overline{O_i} p_{j/i} } -1 \label{eqn:ln_p}\\
   & \Rightarrow p_{m/ij} = \frac{\exp(\theta_m \sum_{q} \beta_q X^q_{ijm})}{\exp(1-\theta_m \frac{\mu_{ij}}{\overline{O_i} p_{j/i}})} \nonumber \\
    & \text{Using that } \sum_{m} p_{m/ij} = 1, \text{we obtain the desirable form.} \\
    & \therefore p_{m/ij} = \frac{\exp(\theta_m \sum_{q} \beta_q X^q_{ijm})}{\sum_{m'} \exp(\theta_{m} \sum_{q} \beta_q X^q_{ijm'})} \label{eqn:p_m} 
\end{align}
\end{subequations}
 
Based on the definition, $S_{ij} = \frac{1}{\theta_m} \ln \sum_{m'} \exp(\theta_{m} \sum_{q} \beta_q X^q_{ijm'})$, we can derive the below equation for the later convenience. 
\begin{align}\label{eqn:ln_p_m}
    \ln p_{m/ij} = \theta_m \sum_{q} \beta_q X^q_{ijm} - \theta_m S_{ij}
\end{align}

Next, from Equation (\ref{eqn:1}), we obtain the following equation by dividing it into $\overline{O_i}$.  

\begin{align*}
    -\frac{1}{\theta_j} [(\ln p_{j/i} + 1)] -\frac{1}{\theta_m} [ \sum_{m}p_{m/ij} \ln p_{m/ij}] +  \sum_{k} \beta_k X^k_{ij} + \sum_{m} p_{m/ij} \sum_{q} \beta_q X^q_{ijm} + \frac{\lambda_i}{\overline{O_i}} = 0
\end{align*}

\noindent Substitute $\ln p_{m/ij}$ with Equation (\ref{eqn:ln_p_m}).

\begin{align*}
    \Rightarrow - \frac{1}{\theta_j} [\ln p_{j/i} + 1]  - \sum_{m} p_{m/ij} \sum_{q} \beta_q X^q_{ijm} + S_{ij}  +  \sum_{k} \beta_k X^k_{ij} + \sum_{m} p_{m/ij} \sum_{q} \beta_q X^q_{ijm} + \frac{\lambda_i}{\overline{O_i}}  = 0
\end{align*}

\begin{subequations}
    \begin{align*}
    \Rightarrow & \ln p_{j/i} = \theta_j [S_{ij} + \sum_{k} \beta_k X^k_{ij} + \frac{\lambda_i}{\overline{O_i}}]  - 1\\
    \Rightarrow & p_{j/i} = \frac{\exp(\theta_j [S_{ij} + \sum_{k} \beta^k X^k_{ij}])}{ \exp(1 - \frac{\lambda_i}{\overline{O_i}})}\\
    & \text{Using that } \sum_{j} p_{j/i} = 1, \text{we obtain the desirable form.} \\
    & \therefore p_{j/i} = \frac{\exp(\theta_j [S_{ij} + \sum_{k} \beta^k X^k_{ij}])}{ \sum_{j'} \exp(\theta_{j} [S_{ij'} + \sum_{k} \beta^k X^k_{ij'}])}
    \end{align*}
\end{subequations}

\end{proof}

\subsection{Proof for hierarhical extended logit model}\label{appendix:hier_extend}

\begin{proof}
    Write the Lagrangian function.

\begin{align*}
    \mathcal{L} =& - \sum_{m} \sum_{a_m} \int^{f^m_a}_{0} g^m_a(w)dw\\
    & - \frac{1}{\theta_r}[\sum_{ijMm} \overline{O_i} p_{j/i} p_{M/ij} p_{m/M} p_{r/ijm} \ln (\frac{p_{r/ijm}}{PS^{ijm}_r})] \\ 
    & - \frac{1}{\theta_j} [\sum_{ij} \overline{O_i} p_{j/i} \ln (p_{j/j}) 
    -\sum_{ij} \overline{T_{ij}} \ln(\frac{\overline{T_{ij}}}{\overline{O_{i}}})] \\
    & - \frac{1}{\theta_m} \left[ \sum_{ijM} \overline{O_i} p_{j/i} p_{M/ij} \ln (p_{M/ij})  
    -\sum_{ijm} \overline{T_{ijM}} \ln(\frac{\overline{T_{ijM}}}{\overline{T_{ij}}}) \right. \\
    & \left. \quad \quad + \sum_{M} \tau_M [\sum_{ijm} \overline{O_i} p_{j/i} p_{M/ij} p_{m/M} \ln (p_{m/M}) 
    -\sum_{ijm} \overline{T_{ijm}} \ln(\frac{\overline{T_{ijm}}}{\overline{T_{ijM}}})] \right] \\
    & + \sum_{k} \beta_k [\sum_{ij} \overline{O_i} p_{j/i} X^k_{ij} - \sum_{ij} \overline{T_{ij}} X^k_{ij}] \\
    & + \sum_{q} \beta_q [\sum_{ijm} \overline{O_i} p_{j/i} p_{M/ij} p_{m/M} X^q_{ijm} - \sum_{ijm} \overline{T_{ijm}} X^q_{ijm}] \\
    & + \sum_{i} \lambda_i [\sum_{j} p_{j/i} - 1 ] \\
    & + \sum_{ij} \mu_{ij} [\sum_{M} p_{M/ij} - 1 ] \\
    & + \sum_{M} \kappa_{M} [\sum_{ijm} p_{m/M} - 1 ] \\
    & + \sum_{ijm} \nu_{ijm} [\sum_{ijmr} p_{r/ijm} - 1 ] 
\end{align*}

Let us focus on the first term. 

\begin{align}
    \frac{\partial \left( \sum_{m} \sum_{a_m} \int^{f^m_a}_{0} g^m_a(w)dw \right)}{\partial p_{r/ijm}} 
    &= \sum_{a_m \ni r} g^m_a(f^m_a) \cdot \frac{\partial f^m_a}{\partial p_{r/ijm}} 
    = \overline{O_i} p_{j/i} p_{M/ij} p_{m/M} \sum_{a_m \in r} g^m_a(f^m_a) \delta^{a_m}_{ijr} \nonumber \\
    & = \overline{O_i} p_{j/i} p_{M/ij} p_{m/M} \cdot g_{ijmr} \\
    \frac{\partial \left( \sum_{m} \sum_{a_m} \int^{f^m_a}_{0} g^m_a(w)dw \right)}{\partial p_{m/M}} 
    & = \sum_{a_m \ni r} g^m_a(f^m_a) \cdot \frac{\partial f^m_a}{\partial p_{m/M}} 
    = \overline{O_i} p_{j/i} p_{M/ij} \sum_{r} p_{r/ijm} \sum_{a_m \in r} g^m_a(f^m_a) \delta^{a_m}_{ijr} \nonumber \\
    & = \overline{O_i} p_{j/i} p_{M/ij} \sum_{r} p_{r/ijm} \cdot g_{ijmr}
\end{align}

\noindent Likewise, 

\begin{align}
    \frac{\partial \left( \sum_{m} \sum_{a_m} \int^{f^m_a}_{0} g^m_a(w)dw \right)}{\partial p_{M/ij}} 
    & = \overline{O_i} p_{j/i} \sum_{m} p_{m/M} \sum_{r} p_{r/ijm} \cdot g_{ijmr}\\
    \frac{\partial \left( \sum_{m} \sum_{a_m} \int^{f^m_a}_{0} g^m_a(w)dw \right)}{\partial p_{j/i}} 
    & = \overline{O_i} \sum_{M} p_{M/ij} \sum_{m} p_{m/M} \sum_{r} p_{r/ijm} \cdot g_{ijmr}
\end{align}

Obtain the first order condition $\frac{\partial \mathcal{L}}{\partial p_{r/ijm}} = 0$. 

\begin{align*}
    \frac{\partial \mathcal{L}}{\partial p_{r/ijm}} &= - \overline{O_i} p_{j/i} p_{m/ij} g_{ijmr}
    - \frac{1}{\theta_r} \overline{O_i} p_{j/i} p_{m/ij} \left( \ln(\frac{p_{r/ijm}}{PS^{ijm}_{r}}) + 1 \right) + \nu_{ijm} = 0 \\
    & \Rightarrow  -g_{ijmr} - \frac{1}{\theta_r} \left[ \ln(p_{r/ijm}) - \ln(PS^{ijm}_r) + 1 \right] + \frac{\nu_{ijm}}{\overline{O_i} p_{j/i} p_{m/ij}} = 0 \\
    & \Rightarrow \ln(p_{r/ijm}) = \theta_r \left( -g_{ijmr} + \frac{\nu_{ijm}}{\overline{O_i} p_{j/i} p_{m/ij}} \right) + \ln(PS^{ijm}_{r}) - 1\\
    & \Rightarrow p_{r/ijm} = \frac{ PS^{ijm}_{r} \exp(-\theta_r g_{ijmr})}{\exp(1 - \frac{\theta_r \nu_{ijm}}{\overline{O_i} p_{j/i} p_{m/ij}})} 
\end{align*}

\noindent Using Equation (\ref{FirstStage:prob4}), we obtain the equation below.
\begin{align}\label{proof_HL:p_rijm}
    \therefore p_{r/ijm} = \frac{ PS^{ijm}_{r} \exp(-\theta_r g_{ijmr})}{\sum_{r} PS^{ijm}_{r} \exp(-\theta_r g_{ijmr})} \quad 
\end{align}

\noindent Using the definition of $S_{ijm} = \frac{1}{\theta_r} \ln \left( \sum_{r} PS^{ijm}_{r} \exp(-\theta_r g_{ijmr}) \right)$ Equation (\ref{proof_HL:p_rijm}) can be written as below for the later convenience.

\begin{align}\label{proof_HL:p_rijm_variant}
    \ln (\frac{p_{r/ijm}}{PS^{ijm}_{r}}) = -\theta_r g_{ijmr} - \theta_r S_{ijm} 
\end{align}

Obtain the first order condition $\frac{\partial \mathcal{L}}{\partial p_{m/M}} = 0$. 

\begin{align*}
    \frac{\partial \mathcal{L}}{\partial p_{m/M}} 
    &=
    - \overline{O_i} p_{j/i} p_{M/ij} \sum_{r} p_{r/ijm} \cdot g_{ijmr} 
    - \frac{1}{\theta_r} \left[ \overline{O_i} p_{j/i} p_{M/ij} \sum_{r} p_{r/ijm} \ln(\frac{p_{r/ijm}}{PS^{ijm}_{r}}) \right] \\
    & - \frac{\tau_M}{\theta_m} \left[ \overline{O_i} p_{j/i} p_{M/ij} \left(\ln(p_{m/M}) + 1\right) \right] + \overline{O_i} p_{j/i} p_{M/ij} \sum_{q} \beta_q X^q_{ijm} + \kappa_M = 0 \\
    & \text{Using Equation (\ref{proof_HL:p_rijm_variant}),} \\
    & \Rightarrow  \frac{\tau_M}{\theta_m} \left(\ln(p_{m/M}) + 1\right) = \sum_{q} \beta_q X^q_{ijm} + \kappa_M + S_{ijm} \\
    & \Rightarrow \ln(p_{m/M}) = \theta_m \left[\sum_{q} \beta_q X^q_{ijm} + S_{ijm}\right]/ \tau_M + \frac{\theta_m \kappa_M}{\tau_M} - 1\\
    & \Rightarrow p_{m/M} = \frac{\exp\left(\theta_m \left[\sum_{q} \beta_q X^q_{ijm} + S_{ijm}\right]/ \tau_M\right)}{\exp(1-\frac{\theta_m \kappa_M}{\tau_M})} 
\end{align*}

\noindent Using Equation (\ref{FirstStage:prob3}),
\begin{align}
    & \therefore p_{m/M} = \frac{\exp \left(\theta_m \left[\sum_{q} \beta_q X^q_{ijm} + S_{ijm}\right]/ \tau_M \right)}
    {\sum_{m'} \exp\left(\theta_m \left[\sum_{q} \beta_q X^q_{ijm'} + S_{ijm'}\right]/ \tau_M\right)} \quad 
\end{align}

\noindent Using the definition of $IV_M = \tau_M \ln \left(\sum_m \exp(\theta_m[\sum_{q} \beta_q X^q_{ijm} + S_{ijm}]/ \tau_M) \right)$, we draw the equation below for later convenience.

\begin{align}
    \ln (p_{m/M}) = \frac{\theta_m \left[\sum_{q} \beta_q X^q_{ijm} + S_{ijm}\right]}{\tau_M} - \frac{IV_{M}}{\tau_M}
\end{align}

Similarly, from $\frac{\partial \mathcal{L}}{\partial p_{M/ij}} = 0$, 
\begin{align}
    \therefore p_{M/ij} = \frac{\exp(IV_M)}{\sum_{M} \exp(IV_M)}
\end{align}

\noindent Using the definition of $S_{ij} = \frac{1}{\theta_m} \ln(\sum_{M} \exp(IV_M))$, we get $\ln p_{M/ij} = IV_M - \theta_m S_{ij}$

Lastly, from $\frac{\partial \mathcal{L}}{\partial p_{j/i}} = 0$, we obtain 

\begin{align}
    \therefore p_{j/i} = \frac{\exp\left(\theta_j (S_{ij} + \sum_{k} \beta_k X^k_{ij})\right)}{\sum_{j} \exp\left(\theta_j (S_{ij} + \sum_{k} \beta_k X^k_{ij})\right)}
\end{align}

\end{proof}








\end{document}